\documentclass{article}

\usepackage{amsfonts,amsmath,graphicx}

\newcommand{\RR}{{\mathbb R}}

\newcommand{\abs}[1]{\lvert#1\rvert}
\newcommand{\norm}[1]{\lVert#1\rVert}

\newcommand{\dv}[2]{{\frac{\partial #1}{\partial #2}}}

\newcommand{\dotex}{{\frac{d}{dt}}}

\newtheorem{thm}{Theorem}
\newtheorem{lem}{Lemma}
\newtheorem{prop}{Proposition}
\begin{document}

\title{Symmetry-based observers for some water-tank problems}

\author{Didier Auroux\thanks{Laboratoire J. A. Dieudonn\'e,
        Universit\'e de Nice Sophia Antipolis, Parc Valrose, 06108 Nice cedex 2, France
        ({\tt auroux@unice.fr}).}
         and Silv\`{e}re Bonnabel\thanks{CAOR, Mines-ParisTech (Ecole des Mines de Paris), 60 Bd St-Michel, 75272 Paris Cedex 06, France ({\tt silvere.bonnabel@mines-paristech.fr}).}}

\markboth{D. Auroux and S. Bonnabel}%
{Symmetry-based observers for some water-tank problems}

\date{}
\maketitle

\begin{abstract}
In this paper we consider a tank containing fluid and we want to estimate the horizontal currents when the fluid surface height is measured. The fluid motion
is described by shallow water equations in two horizontal dimensions.  We build a simple non-linear observer which takes advantage of the symmetries of fluid dynamics laws. As a result its structure is based on convolutions with smooth isotropic kernels, and the observer is remarkably robust to noise. We prove the convergence of the observer around a steady-state. In numerical applications local exponential convergence is expected.  The observer is also applied to the problem of predicting the ocean circulation. Realistic simulations illustrate the relevance of the approach compared with some standard oceanography techniques.
\end{abstract}

\section{Introduction}

The following study is derived from a data assimilation problem in oceanography. The problem considered in this paper consists in estimating the state of a fluid in a water tank where the surface height is measured everywhere. In this paper we propose a symmetry-based non-linear infinite dimensional observer and we prove the convergence  when the fluid motion is described by linearized wave equations under shallow water approximations.

Over the last years much attention has been devoted to the motion planning and feedback stabilization of a fluid under shallow water approximations, problem raised by  \cite{dubois,petit-rouchon-ieee02}. A related problem is the control of flows described by Saint-Venant equations in channels  \cite{coron-ieee,coron-et-al-ecc99,coron-cocv02,mareels,prieur}. Fewer efforts have been put on the theory of observers for this kind of infinite dimensional systems.  Nevertheless a natural extension of this theoretical observer problem consists in oceanographic applications, as we will see later, and extended Kalman filters-type observers are frequently used  to tackle these related problems \cite{Anthes,Verron-Holland,Evensen}. A different approach for observer design for flows in channels is to approximate the motion by non-linear ordinary differential equations at critical points along the channels \cite{mareels,besancon}.  More generally, past efforts in the theory of observers for systems described by partial differential equations (PDEs)  include infinite dimensional Luenberger observers for linear systems \cite{lasiecka,krstic}. Some other problems  have also drawn attention recently \cite{deguenon,deguenon2,Guo,GuoShao}.

Kalman-type filters, or Luenberger observers,  are usually in the standard form ``copy of the system plus injection of the output estimation error (correction term)". For this reason, they
do not take into account the symmetries of the model. There has been recent work on observer design and
symmetries for engineering problems when the model is finite
dimensional and when there is a Lie group acting on the state space
\cite{aghannan-rouchon-ieee03,Aghannan,mahony-et-al-dcd05,arxiv-07,arxiv-08}.
Symmetries provide a helpful guide to design non-linear correction terms. Indeed the only difference
between the observer and model equations comes from the correction
term. Linear systems are invariant by scaling, and so is the correction term in general (Luenberger observer, Kalman filter). But when the system is non-linear, there is no reason why the correction term should have a linear form (extended Kalman filter). When this term  is bound to preserve symmetries, it has a non-linear structure based on the specific nonlinearities of the system, and the observer
is called ``invariant", or ``symmetry-preserving". The result is
that the estimations do not depend on arbitrary choices of units or
coordinates, and the estimates share common physical properties with
the true physical variables (in the examples given in \cite{arxiv-07},
estimated chemical concentrations are automatically positive, estimated
rotation matrices automatically belong to $SO(3)$). In some cases,
the error system even presents very nice properties (autonomous
error equation in \cite{arxiv-08,lageman}).

Looking at \cite{MarSal2007a,bonnabel_auto,arxiv-07}, the design method of symmetry-preserving observers could be summed up this way: the non-linear form of the observer is given by the symmetries, and the gains are tuned assigning the poles of the error system around a trajectory or a steady-state. This is always possible as around any steady-state, invariant observers can be identified to  Luenberger observers \cite{arxiv-07}.

This paper is an extension to the infinite-dimensional case of the
recent ideas on observer design and symmetries for systems described
by ordinary differential equations (ODEs). The Saint-Venant equations considered in this paper are indeed invariant by rotation and translation
$(SE(2)$-invariance). In the case of systems described by PDEs, the design of
observers based on the symmetries of the physical system is new to
the authors' knowledge.

The first theoretical contribution of this paper is to derive a $SE(2)$-invariant observer for the problem. The correction terms do not depend on any
non-trivial choice of coordinates. They correspond to a convolution product of the output error and a smooth isotropic kernel, a feature which ensures remarkable robustness to white noise.
  With respect to this latter feature, the observer is close in its flavour to \cite{krstic-05} where the authors derive a non-linear observer to estimate the velocity and pressure in an infinite channel. Their observer consists in a copy of the system and a correction term corresponding to a one-dimensional convolution product of the output error and some kernel (but in one dimension there is no such thing as invariance by rotation).
We then yield proof of
convergence of the estimation error to zero for the linearized Saint-Venant equations.  The form of the error equation shows that in numerical experimentations exponential convergence can always be expected. A noticeable fact is that the observer
depends only on a small number of parameters, as the respect of the symmetries implies some restrictions on the non-linear correction terms. All these
parameters admit a physical interpretation. So the observer gains
are a lot easier to tune than the Kalman gains, and the corresponding computing cost is very low compared to the extended Kalman filter.

 The second theoretical contribution is to extend those results to a large class of $SE(2)$-invariant observers with smoothing correction terms.  The idea to derive
systematical smoothing terms based on physical symmetries, is
standard in image processing, and was initiated by \cite{lions}. One
of the simplest  method consists indeed in a convolution with a two-dimensional smooth
rotation-invariant kernel (isotropic diffusion based on the heat
equation), but we prove convergence for a much larger class of rotation-invariant kernels.  The major difference with  \cite{lions} is that those smoothing terms are combined with a
dynamical model to provide an estimation of physical quantities
which are not directly measured, i.e. we build  observers.

The observer is applied to an oceanography example.
The problem considered is the following: the ocean is described by a
simplified shallow-water model \cite{jiang}. The sea surface height (SSH) is
measured (with noise) by satellites. The goal is to
estimate the height, and the marine currents (not measured).
Observers are a flexible assimilation technique,
computationally much more economical than variational data
assimilation methods \cite{LeDimet,LeDimetTalagrand}. Nevertheless as the oceanographic models have become very complex in
the recent years, the high computing cost of most extended Kalman
filters (EKF) can still be prohibitive for data assimilation \cite{Rozier07}. Our observer is thus a relevant challenger to the EKF, due to its structural properties (geometric structure, local convergence proof, isotropic smoothing term, easiness of the gain tuning), and its low computational cost, as illustrated by extensive numerical experiments.

The paper is organized as follows. In Section \ref{saint-venant:sec} we define a class of $SE(2)$-invariant
observers for water-tank systems for which the surface height is measured.
Convergence of the estimation error to zero on the first-order
approximation of this system is proved. In Section \ref{oceanography:sec} we consider a bi-dimensional
shallow water model, often used in geophysics for ocean or fluid flow modeling. We propose a $SE(2)$-invariant non-linear observer when the SSH is measured by satellites. In Section \ref{num:sec}, we report the results of extensive
numerical simulations on both the linearized  and
nonlinear shallow water models, which illustrate the properties of the observer. These results are compared with another
standard oceanography technique based on the use of observers (nudging). We also show the
remarkable robustness of the estimation to Gaussian white noise on the observations.
Finally, some conclusions and perspectives are given in Section \ref{ccl:sec}.

\section{Water-tank system, symmetries, and observer design}\label{saint-venant:sec}

The problem we are concerned with is the motion of a perfect fluid under gravity described by Saint-Venant equations with a free surface (the shallow water assumption). The state of
the fluid is its surface height, and the horizontal speed of the
currents. The choice of the orientation  and the origin of the frame
of $\RR^2$  used to express the horizontal coordinates
$(x,y)\in\RR^2$ is arbitrary: the physical problem is invariant by
rotation and translation. Indeed from a mathematical viewpoint
the Laplace operator $\Delta$ is invariant by rotation and
translation.  The first term of any observer for this problem is
automatically invariant by rotation and translation, as it is a copy
of the equations of the physical system. There is no reason why the
correction term should depend on any non-trivial choice of the
orientation and origin of the frame. It would yield correction terms
giving more importance to the values of the height measured in some
arbitrary direction of $\RR^2$. In the general case, without
additional information on the model, it seems perfectly logical
 to correct  the observer isotropically. This constraint
suggests interesting correction terms.

\subsection{Saint-Venant model}\label{sec:SV}

Consider a rectangular domain: $0\le x\le L$
and $0\le y\le L$ (which can be considered square without loss of generality), where $x$ and $y$ are cartesian coordinates. Let $\nabla$ be the
corresponding gradient operator:
$$
\nabla=\left(\frac{\partial}{\partial x}, \frac{\partial}{\partial
y}\right)^T.
$$
The Saint-Venant equations write:
\begin{eqnarray}
\dv{h}{t}&=&-\nabla\cdot(hv), \label{saint-venant:eq1}\\
\dv{v}{t}&=&-(v\cdot\nabla) v-g\nabla h. \label{saint-venant:eq2}
\end{eqnarray}
where $hv=h(v_x \textbf{i} + v_y \textbf{j})$ is the horizontal
transport, with $\textbf{i}$ and $\textbf{j}$ denoting the axes of an Euclidean frame  and $g$ is the gravity.

The boundary conditions that we consider are:
\begin{itemize}
\item rigid boundaries: $v_x(x,y)=0$ for $x=0$ and $x=L$, $\forall y$; and $v_y(x,y)=0$ for $y=0$ and $y=L$, $\forall x$. In other words, $v.\textbf{n}=0$ on the boundary of the domain, $\textbf{n}$ being the outward unit normal to the domain.
\item no-slip lateral boundary conditions for $v_x$ on the top and bottom boundaries of the domain, and for $v_y$ on the left and right boundaries. As the domain does not move, the no-slip lateral conditions are equivalent to $v_x(x,y)=0$ for $y=0$ and $y=L$, $\forall x$; and $v_y(x,y)=0$ for $x=0$ and $x=L$, $\forall y$.
\end{itemize}
All together, the boundary conditions are $v=0$. The theory of characteristics in 2D tells us that in our case (no normal velocity through the boundary), only two boundary conditions must be imposed on each boundary. Then, there is no need for boundary conditions on the height, as equation \eqref{saint-venant:eq1} is a standard transport equation. The initial conditions $(h(0),v(0))$ complete the system.

Note that from a computational point of view, the equations are discretized on an Arakawa C grid \cite{arakawa}, in which the velocity components are defined at the center of the edges. Then, instead of imposing $v_x=0$ e.g. on the top boundary ($y=L$), a classical way to impose no-slip boundary conditions is to use an additional row of points in the grid beyond the boundary, on which $v_x$ is set to the opposite of the value of $v_x$ on the first inside row of points, so as to ensure a null mean value at the boundary. Another standard set of boundary conditions is rigid boundaries and free-slip lateral boundary conditions \cite{adcroft}.

We assume that the height $h(x,y,t)$ is measured (with noise) for all $x,y,t$.
The problem is  the estimation of $v(x,y,t)$ at
any point $(x,y)\in [0,L]^2$ of the domain. In the presence of noise, the problem is the estimation of both variables $v$ and $h$.

Note that this assumption could be slightly relaxed. In oceanographic applications  discussed in Section \ref{oceanography:sec}, the height can be partially observed. But then, after some time, all the observations are usually gathered together on a same observation map, and then interpolated in order to obtain full spatial observations of $h$ (for all $x$ and $y$), at some discrete times $t$. We could then assume that $h(x,y,t)$ is known for all $x$ and $y$, but at only some times $t$. Such an approach has the main advantage of allowing one to spatially filter the data, and thus allows the method proposed in this paper to be applied.  In the case of discrete (in space) sets of measurements, the  correction terms proposed below, with no measurements at some locations,  boil down to standard Luenberger observer-like correction terms.

\subsection{Model symmetries}

The unit vectors $\textbf{i}$ and $\textbf{j}$ can be chosen to point East and
North respectively. This choice is arbitrary, and the equations of
fluid mechanics depend neither on the orientation nor on the
origin of the frame in which the coordinates are expressed: they are
invariant under the action of the Lie group $SE(2)$, the {\em Special Euclidean}
group of isometries of the plane $\mathbb{R}^2$. Let us prove it. Let $R_\theta=
 \left(
    \begin{array}{cc}
      \cos\theta & -\sin\theta \\
      \sin\theta & \cos\theta \\
    \end{array}
  \right)
 $
be a horizontal rotation of angle $\theta$. Let
$(x_0,y_0)\in\mathbb{R}^2$ be the origin of some new frame. Let
$(X,Y)$ be the coordinates associated to this new frame
$R_{-\theta}(\textbf{i},\textbf{j})-(x_0,y_0)$. In this new frame,
the variables read
\begin{eqnarray}
(X,Y)&=&R_\theta(x,y)+(x_0,y_0),\label{new:coord:eq1}\\
H(X,Y,t)&=&h(x,y,t),\label{new:coord:eq2}\\
V(X,Y,t)&=&R_\theta v(x,y,t).\label{new:coord:eq3}
\end{eqnarray}
and $(\dv{}{X},\dv{}{Y})=R_\theta(\dv{}{x},\dv{}{y})$ which implies
$\nabla H (X,Y,t)=R_\theta \nabla h (x,y,t)$. The equations in the new
coordinates are unchanged:
\begin{eqnarray}
\frac{\partial H}{\partial t}&=&-\nabla\cdot (HV).\\
\dv{V}{t}&=&-(V\cdot\nabla) V-g\nabla H.
\end{eqnarray}
The Laplace and divergence operators are unchanged by the
transformation as they are invariant to rotations (although they are
usually written in fixed coordinates, their value do not depend on
the orientation of the chosen frame). Note that the square domain
$D=[0,L]^2\subset\RR^2$ has to be replaced by the square
$\left(R_\theta D+(x_0,y_0)\right)\subset\RR^2$. The boundary conditions ($v=0$ on the boundaries of the domain) are also clearly invariant by rotation.

\subsection{A  symmetry-preserving observer}\label{obs:sec}

Any non-linear observer for the system \eqref{saint-venant:eq1}-\eqref{saint-venant:eq2}
  writes :
\begin{eqnarray}
\frac{\partial \hat h}{\partial t}&=&-\nabla\cdot (\hat h\hat
v)+F_h(h,\hat v,\hat h), \label{obs0:eq1}\\
\dv{\hat v}{t}&=&-(\hat v\cdot\nabla) \hat v-g\nabla \hat h+F_{v}(h,\hat v,\hat h), \label{obs0:eq2}
\end{eqnarray}
with the boundary condition $\hat v=0$ on all boundaries of the domain, and where the
correction terms vanish when the estimated height $\hat h$ is equal
to the observed height $h$:
$$F_{v}(h,\hat v,h)=0,\quad F_{h}(h,\hat v,h)=0.$$
We propose the following observer for the system
\eqref{saint-venant:eq1}-\eqref{saint-venant:eq2}:
\begin{eqnarray}
\frac{\partial\hat h}{\partial t}&=&-\nabla\cdot(\hat h\hat v) \nonumber \\
&&+\iint\phi_h(\xi^2+\zeta^2)\ (h-\hat h)_{(x-\xi,y-\zeta,t)}~d\xi d\zeta \nonumber \\
 &=& -\nabla\cdot(\hat h\hat v)+ \varphi_h * (h-\hat h), \label{observateur:eq1} \\
\frac{\partial \hat v}{\partial t}&=&-(\hat v\cdot\nabla) \hat v-g\nabla
\hat h \nonumber \\
&&+  \iint\phi_{v}(\xi^2+\zeta^2)\ \nabla(h-\hat h)_{(x-\xi,y-\zeta,t)}~d\xi d\zeta \nonumber \\
 &=& -(\hat v\cdot\nabla) \hat v-g\nabla \hat h+\varphi_{v} * \nabla (h-\hat h), \label{observateur:eq2}
\end{eqnarray}
with the same boundary conditions as before, and where
\begin{eqnarray}
\varphi_v(x,y) &=& \beta_v \exp(-\alpha_v (x^2+y^2)), \label{gains:eq1} \\
\varphi_h(x,y) &=& \beta_h \exp(-\alpha_h (x^2+y^2)),
\label{gains:eq2}
\end{eqnarray}
Such an observer preserves the symmetries of the system  as the correction terms are based on a convolution product (translation invariance) with an isotropic kernel (rotation invariance). This is a very logical choice. Indeed,  why should the quality of the estimation
depend upon any non-trivial choice of orientation and origin of the frame when
the physical system under consideration does not at all?

Such correction terms make the observer very robust to noise, as they operate a smoothing of the measured image. The high frequencies in the signal are thus
efficiently filtered.
 Indeed, such translation and rotation invariant terms are standard for image smoothing (see, e.g., \cite{lions}). Other symmetry-preserving smoothing terms will be found in subsection \ref{class:sec}.

While approaching the boundary of the domain, the integrals and convolution kernels may become undefined. But they can easily be extended close to the boundary by truncating the integrals so that they only cover the domain, or equivalently by extending the functions by $0$ outside the domain.

\subsection{Convergence study on the linearized system}\label{cv:sec}

As it seems out of reach to study the convergence of the full
non-linear system, we are going to linearize the system
\eqref{saint-venant:eq1}-\eqref{saint-venant:eq2}  around the steady-state $h=\bar h$ and $v=\bar
v$, using exactly the same simplifications as
\cite{petit-rouchon-ieee02} which considers the
open-loop control problem of system
\eqref{saint-venant:eq1}-\eqref{saint-venant:eq2} with boundary
control. The considered equilibrium is characterized by $\bar h$
equal to a constant height, and $\bar v = 0$. The observer gains are
designed on this latter system, and we prove at the end of this
section that they ensure the strong asymptotic convergence of the
error.

Approximating the true system with  the linearized system means that we only
consider small velocities $\delta v=v -\bar v \ll \sqrt{g\bar h}$
and heights $\delta h=h-\bar h \ll \bar h$. Note that these
approximations are consistent with the first set of numerical experiments
(subsection \ref{num:linearized:sec}), in which the ratio $\delta v$ (resp. $\delta
h$) to $\sqrt{g\bar h}$ (resp. $\bar h$) is of the order of
$10^{-2}$ to $10^{-3}$. The linearized system is
\begin{eqnarray}
\dv{(\delta h)}{t}&=&-\bar h\ \nabla \cdot\delta v, \\
\dv{(\delta v)}{t}&=&-g\nabla \delta h,
\end{eqnarray}
and the estimation errors, $\tilde h = \delta \hat h - \delta h $
and $\tilde v = \delta \hat v - \delta v$, are solution of the
following linear equations:
\begin{eqnarray}
\dv{\tilde h}{t}&=&-\bar h\ \nabla\cdot \tilde v-\varphi_h * \tilde h,\\
\dv{\tilde v}{t}&=&-g\nabla \tilde h -\varphi_v*\nabla \tilde h.
\end{eqnarray}
Eliminating $\tilde v$ and using $\nabla(\varphi_v*\nabla
h)=\varphi_v*\Delta h$  yields a modified damped wave equation with
external viscous damping:
\begin{equation}
\label{wav:eq} \frac{\partial^2 \tilde h}{\partial t^2}= g\bar
h\Delta \tilde h + \bar h\, \varphi_v*\Delta \tilde h - \varphi_h*
\frac{\partial \tilde h}{\partial t}.
\end{equation}
\begin{thm}
\label{thm1} If $\varphi_v$ and $\varphi_h$ are defined by
(\ref{gains:eq1}) and (\ref{gains:eq2}) respectively with
$\beta_v,\beta_h,\alpha_v,\alpha_h>0$, then the first order
approximation of the error system around the equilibrium
$(h,v)=(\bar h,0)$ given by \eqref{wav:eq} is strongly
asymptotically convergent. Indeed if we consider the following
Hilbert space and norm: $\mathcal{H}=H^1(\Omega)\times L^2(\Omega)$,
\begin{equation}
\norm{(u,w)}_\mathcal{H}= \left(\int_\Omega \norm{\nabla
u}^2+\abs{w}^2\right)^{1/2},
\end{equation}
then, for every $\tilde{h}$ solution of \eqref{wav:eq},
\begin{equation}\label{eq:thm}
\lim_{t\rightarrow\infty}\ \bigg\|\left(\tilde
h(t),\frac{\partial\tilde h}{\partial
t}(t)\right)\bigg\|_\mathcal{H}=0\,.
\end{equation}
\end{thm}
This theorem proves the strong and asymptotic convergence of the
error $\tilde{h}$ towards $0$, and then it also gives the same
convergence for $\tilde{v}$. We deduce that the observer
\eqref{observateur:eq1}-\eqref{observateur:eq2} tends to the true
state when time goes to infinity.

A dimensional analysis can yield a meaningful choice of the gains.
The parameters $\alpha_v^{-2}$, $\alpha_h^{-2}$ are expressed in
meters. They define the size of the regions of influence of the
kernels, i.e. the region around any point in which the measured
values of $h$ are used to correct the estimation at the point. These
values can be set experimentally using the data from the physical
system. Moreover,  $\beta_v$ and $\beta_h$
can be tuned via  the following heuristics. The error system
\eqref{wav:eq} can be approximated by the following system, which
corresponds to the case $\alpha=+\infty$:
\begin{equation}
\label{wav2:eq} \frac{\partial^2 \tilde h}{\partial
t^2}+2\xi_0\omega_0 \frac{\partial \tilde h}{\partial t}=
(L_0\omega_0)^2\Delta \tilde h .
\end{equation}where $L_0^2\omega_0^2=g\bar h+\bar h\beta_v$,
$2\xi_0\omega_0=\beta_h$, as long as we impose $L_0^2\omega_0^2\geq
g\bar h$. $\beta_v$ and $\beta_h$ can be chosen in order to control
the characteristic pulsation
  $\omega_0$, length $L_0$, and damping coefficient $\xi_0$ of the approximated
   error equation \eqref{wav2:eq}. These quantities have an obvious physical meaning and
    can be set accordingly to the characteristics of  the physical system under consideration.
     Such heuristics provide a first reasonable tuning of the gains.

\subsection{Proof of theorem \ref{thm1}} \label{sec:proof}

In this section,  we prove the
strong convergence of the error system in the Hilbert space
$\mathcal{H}$. The proof is inspired by \cite{komornik-book-2} (see also \cite{krstic-05} on an infinite 1D domain). Let $\psi_v=g\bar h\delta_0+\bar h\varphi_v$. For
simplicity reasons, we assume that $L=\pi$. The error equation
\eqref{wav:eq} can be rewritten as a modified wave equation on a
square domain with Dirichlet boundary condition:
\begin{equation}\label{error:eq}
\begin{array}{rclp{0.2cm}l}
\displaystyle \frac{\partial^2 }{\partial t^2}u & \!\!\!=\!\!\! & \psi_v*\Delta
u-\varphi_h*\displaystyle \frac{\partial}{\partial
t}u & & \text{in}\ \ \RR^+\times\Omega,\\
u & \!\!\!=\!\!\! & 0 & & \text{on}\ \,\RR^+\times\partial\Omega,\\
u(0) & \!\!\!=\!\!\! & u_0, \ \  u_t(0)=u_1 & & \text{in}\ \ \Omega,
\end{array}
\end{equation}
where $\Omega=[0,\pi]^2$, and $u(t,x,y)$ represents the estimation error $\tilde h$.

We denote by $(e_{pq})$ the following orthonormal basis of
$H^1_0(\Omega)$, composed of eigenfunctions of the unbounded
operator $\Delta$:
\begin{equation}\label{eq:epq}
e_{pq}=\frac{2}{\pi}\sin(px)\sin(qy).
\end{equation}
Moreover, let $f(s)=(2\beta_v)^{1/4} \exp(-2\alpha_vs^2)$ and  $g(s)=(2\beta_h )^{1/4}\exp(-2\alpha_h s^2)$. As the convolution product of two Gaussians is a Gaussian we have
\begin{eqnarray}
\varphi_v(x,y)&=&(f(x)*f(x))(f(y)*f(y)), \label{gains:ff}\\
\varphi_h(x,y)&=&(g(x)*g(x))(g(y)*g(y)), \label{gains:gg}
\end{eqnarray}
As $f$ and $g$ are even functions, their Fourier
coefficients are real. If we denote by $(\hat f_p)$ and $(\hat g_p)$
the Fourier coefficients of $f$ and $g$ respectively, then, as the convolution is a multiplication in the frequency domain, the
Fourier coefficients of $\psi_v$ are $g\bar h + \bar h \hat f_p^2
\hat f^2_q$. Similarly, the Fourier coefficients of $\varphi_h$ are
$\hat g_p^2 \hat g_q^2$. As all these coefficients are real and
positive, we denote them by $f^2_{pq}$ for $\psi_v$, and $g^2_{pq}$
for $\varphi_h$.  We now need the following
intermediate result:
\begin{lem}\label{lem1}
If $u_0\in H^1_0(\Omega)$ and $u_1\in L^2(\Omega)$, then equation
\eqref{error:eq} has a unique solution satisfying
\begin{equation}
u\in C(\RR^+;H^1_0(\Omega))\cap C^1(\RR^+;L^2(\Omega)).
\end{equation}
It is given by the series:
\begin{equation}\label{eq:upq0}
u(t,x,y)=\frac{2}{\pi}\sum_{1\leq p,q}u_{pq}(t)\sin(px)\sin(qy),
\end{equation}
where $u_{pq}$ can be written either in the following way
\begin{equation}\label{eq:upq1}
u_{pq}(t)=e^{\frac{-g_{pq}^2}{2}t}(A_{pq}\cos(\omega_{pq}t)+B_{pq}\sin(\omega_{pq}t)),
\end{equation}
or
\begin{equation}\label{eq:upq2}
u_{pq}(t)=e^{\frac{-g_{pq}^2}{2}t}(A_{pq}\cosh(\tilde\omega_{pq}t)+B_{pq}\sinh(\tilde\omega_{pq}t)).
\end{equation}
Moreover, the latter case appears at most for a finite number of
indices, and $\tilde\omega_{pq}<\frac{g_{pq}^2}{2}$ (we refer to equation (\ref{eq34}) for the expression of $\omega_{pq}$ and $\tilde{\omega}_{pq}$).
\end{lem}

The proof of the Lemma is as follows. We rewrite equation
\eqref{error:eq} as
\begin{equation}\label{eq:AU}
{\dotex U=\mathcal A U},
\end{equation}
where $U=(u,u_t)$ and $\mathcal A$ is the following unbounded linear
operator on $\mathcal{H}$:
\begin{equation}\label{eq:A}
\mathcal{A}(u,w):=(w, \psi_v*\Delta u-\varphi_h*w).
\end{equation}
From \eqref{eq:A} and \eqref{eq:epq}, we deduce that
\begin{equation}
E_{pq}=\begin{pmatrix}1\\\lambda_{\pm pq}\end{pmatrix}e_{pq}
\end{equation}
are eigenvectors of $\mathcal{A}$ associated to the eigenvalues
$\lambda_{\pm pq}$, solutions of
\begin{equation}\label{eq:lambda}
\lambda_{\pm pq}^2+g_{pq}^2\lambda_{\pm pq}+f_{pq}^2(p^2+q^2)=0.
\end{equation}
Moreover, the family of eigenvectors $(E_{pq})$ forms a Riesz basis
of the Hilbert space $\mathcal{H}$. The discriminant of
\eqref{eq:lambda} is $\Delta_{pq}=g_{pq}^4-4(p^2+q^2)f_{pq}^2$. It
can be positive for a finite number of indices only, since
$g_{pq}^2\rightarrow 0$ and $f_{pq}^2\geq g\bar h$ when $p$ and $q$
go to infinity. We found a Riesz basis of $\mathcal H$ formed by
eigenvectors of $\mathcal A$, the eigenvalues have no finite
accumulation point and their real part are bounded. Thus all
assumptions of theorem 3.1 of \cite{komornik-book-2} are satisfied:
the solution $U$ of \eqref{eq:AU} is given by the series
\begin{eqnarray}
&&\hspace*{-0.5cm}U(t)=\hspace*{-0.3cm}\sum_{\begin{array}{c}\scriptstyle p,q\ge 1\\\scriptstyle\Delta_{pq}<0\end{array}} \hspace*{-0.2cm}\left( U_{pq}e^{\frac{-g_{pq}^2+i\sqrt{4(p^2+q^2)f_{pq}^2-g_{pq}^4}}{2}t} \right. \nonumber \\[-1cm]
&& \hspace*{2.5cm}\left. +U_{-pq}e^{\frac{-g_{pq}^2-i\sqrt{4(p^2+q^2)f_{pq}^2-g_{pq}^4}}{2}t}\right) E_{pq} \nonumber \\
&&\hspace*{-0.3cm}+\hspace*{-0.3cm}\sum_{\begin{array}{c}\scriptstyle
p,q\ge 1\\\scriptstyle\Delta_{pq}\ge0\end{array}}
\hspace*{-0.2cm}\left(
U_{pq}e^{\frac{-g_{pq}^2+\sqrt{g_{pq}^4-4(p^2+q^2)f_{pq}^2}}{2}t}
\right. \nonumber \\[-1cm]
&& \hspace*{1.7cm}\left. +U_{-pq}e^{\frac{-g_{pq}^2-\sqrt{g_{pq}^4-4(p^2+q^2)f_{pq}^2}}{2}t}\right)
E_{pq}. \label{eq34}
\end{eqnarray}
Finally, the coefficients can be found using the Fourier series of
the initial condition. We have
\begin{eqnarray}
A_{pq}&=&\frac{4}{\pi^2} \int_{[0,\pi]^2} u(0,x,y) \sin(px) \sin(qy)\,dxdy,\\
B_{pq}&=&\frac{4}{\omega_{pq}\pi^2}\int_{[0,\pi]^2}
\left(u_t(0,x,y)+\frac{g_{pq}^2}{2}u(0,x,y)\right) \nonumber \\
&& \hspace*{2.5cm}\times \sin(px)\sin(qy)\,dxdy.
\end{eqnarray}

Now the Lemma is proved. All we have to prove now is that the
solution, which represents the estimation error, converges to 0 when
time goes to infinity. Recall that the coefficients $u_{pq}$ are
given by equation \eqref{eq:upq1}, except for a finite number of
indices. Define
\begin{eqnarray}
u_N(t,x,y)&=&\frac{2}{\pi}\sum_{p+q\geq
N}e^{\frac{-g_{pq}^2}{2}t}(A_{pq}\cos(\omega_{pq}t) \nonumber \\
&&\hspace*{0.3cm}+B_{pq}\sin(\omega_{pq}t))\sin(px)\sin(qy).
\end{eqnarray}
Since $u_0\in H^1_0(\Omega)$ and $u_1\in L^2(\Omega)$, Parseval's
theorem tells us that for any $\varepsilon>0$, there exists $N$ such
that
\begin{equation}
\bigg\|u_N(t),\frac{\partial u_N}{\partial t}(t)\bigg\|_{\mathcal
H}\leq\varepsilon /2,\quad \forall t\ge 0.
\end{equation}
From \eqref{eq:upq1} and \eqref{eq:upq2}, there exists $T>0$ such
that for any $t\ge T$,
\begin{equation}
\bigg\|(u-u_N)(t),\frac{\partial (u-u_N)}{\partial
t}(t)\bigg\|_{\mathcal H}\leq\varepsilon /2.
\end{equation}
Finally, $\norm{u,u_t}_{\mathcal H}<\varepsilon$ for any $t\ge T$.
We proved equation (\ref{eq:thm}), i.e. the strong convergence of
the linearized error system.

Note that this proves the result for any kernel functions $\varphi_h$ and $\varphi_v$  provided they are smooth, and their Fourier coefficients be real and strictly positive. Note also that for
$N>0$ arbitrary large, from Lemma \ref{lem1}, the truncated solution
$u_N$ tends to $0$ exponentially in time. Thus exponential
convergence is expected in numerical experiments.

\subsection{A class of locally converging symmetry-preserving observers}\label{class:sec}

This subsection can be skipped by the uninterested reader. Observer \eqref{observateur:eq1}-\eqref{observateur:eq2} preserves the symmetries of the system, it is robust to noise, and it is such that the linearized error equation around fluid at rest converges to zero. However there are many other observers having those desirable properties. In the seminal paper \cite{lions}, the authors seek image-processing transforms that satisfy a list of formal requirements such as translation and rotation invariance.  Inspiring from this work and also from \cite{arxiv-07}, we are going to seek a  class of non-linear observers (i.e. correction terms)
 that   satisfy the following list of formal requirements:
 \begin{itemize}
 \item ``symmetry preservation requirement": invariance to translations and rotations.
 \item ``smoothing by convolution requirement":   to reduce  the noise, the measured  output must be smoothed (especially before being differentiated).
   \item ``local stability requirement": strong asymptotic convergence of the linearized error system.
 \end{itemize}This classification yields a new class of candidate observers which are sensible alternatives to observer \eqref{observateur:eq1}-\eqref{observateur:eq2}. Indeed consider for instance:
\begin{eqnarray}\label{obs:delta:eq}\frac{\partial\hat h}{\partial t}&=& -\nabla\cdot(\hat h\hat v)- \varphi_h * \Delta (h-\hat h), \\
\frac{\partial \hat v}{\partial t}
 &=& -(\hat v\cdot\nabla) \hat v-g\nabla \hat h+\varphi_{v} * \nabla (h-\hat h)
\label{obs:delta:eq2}
\end{eqnarray}
  Such a structural damping term changes drastically the spectrum, and the differentiation process of the measured signal is carried out without amplifying high frequencies (noise). The linearized error equation is then $
\frac{\partial^2}{\partial t^2}\tilde h= (g\bar
h\delta_0+\varphi_v)*\Delta \tilde h +\varphi_h * \Delta
\left( \frac{\partial}{\partial t}\tilde h  \right)$, so $\tilde h$ is given by the series \eqref{eq:upq0} along with \eqref{eq:upq1}-\eqref{eq:upq2} where   $g_{pq}^2$ is replaced everywhere by $g_{pq}^2(p^2+q^2)$. Thus the convergence rate is speeded up by a factor $p^2+q^2$ on each Fourier coefficient, and the high frequencies are still efficiently filtered as the correction terms are automatically smooth. Moreover,  the quality of the estimation does not
depend upon any non-trivial choice of orientation and origin of the frame.  Indeed, the Laplace operator is $SE(2)$-invariant.

\paragraph{Symmetry preservation}In fact,
 according to standard results (see e.g. \cite{schwartz-seminaire}), \emph{any}  $SE(2)$-invariant scalar differential  operator
 writes $Q(\Delta)$, where $Q$
is a polynomial and $\Delta$ is the Laplacian.
 To fill the first requirement, this feature suggests  to use polynomials of the Laplacian to design correction terms for the general form \eqref{obs0:eq1}-\eqref{obs0:eq2}. To get a symmetry-preserving scalar correction term  $F_{h}(h,\hat v,\hat h)$, the coefficients of the polynomials must depend on invariant scalar functions of $h,\hat v,\hat h$. Thus they must depend on $\hat v$ only via an invariant function of
$\hat v$, typically $|\hat v|^2$.  A large class of symmetry-preserving correction terms is:
\begin{equation}
F_h = Q_1(\Delta,h,|\hat v|^2,\hat h-h) + \nabla
\left(Q_2(\Delta,h,|\hat v|^2,\hat h-h)\right)\cdot \hat v,
\end{equation}
where $Q_1$ and $Q_2$ are scalar polynomials in $\Delta$. More precisely, for $i=1,2$, we
have
\begin{eqnarray}
\hspace*{-0.8cm}Q_i(\Delta,h,|\hat v|^2,\hat h-h)\!\!&=&\!\! \sum_{k=0}^{N} a^i_k(h,|\hat
v|^2,\hat h-h) \nonumber\\
&&~\Delta^{k} \left(b^i_k(h,|\hat v|^2,\hat
h-h)\right),
\end{eqnarray}
where $a_k^i$ and $b_k^i$ are smooth scalar functions such that $
a_k^i(h,|\hat v|^2,0)=b_k^i(h,|\hat v|^2,0)=0.
$
For the vectorial correction term $F_{v}$, we use the vectorial
counterpart of $F_h$:
\begin{equation}
F_{v}= P_1(\Delta,h,|\hat v|^2,\hat h-h)~ \hat v +  \nabla
\left(P_2(\Delta,h,|\hat v|^2,\hat h-h) \right),
\end{equation}
where $P_1$ and $P_2$ are polynomials in $\Delta$, like $Q_1$ and
$Q_2$.
\paragraph{Symmetry preservation and smoothing by convolution} The polynomials above involve a differentiation process, and thus must be coupled with a filtering process. Let us  find integral terms $F_h$ and $F_{v}$ that are $SE(2)$-
invariant. They can be expressed as a
convolution between the previous invariant differential terms, and a
two-dimensional kernel $\psi(\xi,\zeta)$. As the correction terms above are invariant to rotation, the value of the kernel should not depend on
any particular direction either, so $\psi$ must be a function of the
invariant $\xi^2+\zeta^2$ (isotropic gain). If we let $\phi_v$
and $\phi_h$ be two real-valued kernels, a class of symmetry-preserving integral correction
terms is:
\begin{eqnarray}
F_{v}(x,y,t) &=& \iint \phi_{v}(\xi^2+\zeta^2) \left[R_1(\Delta,h,|\hat v|^2,\hat h-h) \hat v \right. \nonumber \\
&&  \hspace*{-2cm}+ \left. \nabla \left(R_2(\Delta,h,|\hat v|^2,\hat h-h)\right)\right]_{(x-\xi,y-\zeta,t)} ~  d\xi d\zeta, \label{fv:eq}\\
F_{h}(x,y,t) &=& \iint \phi_{h}(\xi^2+\zeta^2) \left[S_1(\Delta,h,|\hat v|^2,\hat h-h) \right. \nonumber \\
&&  \hspace*{-2cm}+ \left. \nabla \left(S_2(\Delta,h,|\hat v|^2,\hat
h-h)\right)\cdot\hat v\right]_{(x-\xi,y-\zeta,t)} ~  d\xi d\zeta,\label{fh:eq}
\end{eqnarray}
where the polynomials $R_i$ and $S_i$ are defined like the $Q_i$'s.

The support of $\phi_{v}$ (resp. $\phi_{h}$) is a subset of $\RR$.
Its characteristic size defines a zone in which it is significant to
correct the estimation with the measurements. The  observer is independent of any arbitrary
choice of orientation (rotation invariance), as well as of the origin of the chosen
frame (translation invariance). If the kernels are smooth, the correction terms are automatically smooth even if the
measurements are not (noise robustness). Note that, if $\phi_v$ and $\phi_h$ are set equal to
Dirac functions, one recovers the differential terms above.
\paragraph{Local convergence} Although the stability analysis of all symmetry-preserving observers with general correction terms \eqref{fv:eq}-\eqref{fh:eq} is out of reach,  the following proposition, applying to observers \eqref{observateur:eq1}-\eqref{observateur:eq2}-\eqref{gains:eq1}-\eqref{gains:eq2}, and \eqref{obs:delta:eq}-\eqref{obs:delta:eq2},   generalizes Theorem \ref{thm1}   to a large class of observers:
\begin{prop}Let f, g be any smooth functions, and $\varphi_h$ and $\varphi_v$ be given by \eqref{gains:ff}-\eqref{gains:gg}. For any integer $N\geq 0$, if $\lambda_k$ are positive numbers for $0\leq k\leq N$, and at least one of them is strictly positive, the following  class of observers is such that the first order approximation of the error system around $(\bar h,0)$ is strongly asymptotically convergent:
\begin{eqnarray}
\frac{\partial\hat h}{\partial t}&=& -\nabla\cdot(\hat h\hat v)+ \varphi_h *\bigl( \sum_{k=0}^N(-1)^k\lambda_k\Delta^k (h-\hat h)\bigr), \nonumber\\
\frac{\partial \hat v}{\partial t}
 &=& -(\hat v\cdot\nabla) \hat v-g\nabla \hat h+\varphi_{v} * \nabla (h-\hat h) \nonumber
\end{eqnarray}
Moreover, if $\varphi_h(x,y)$ and $\varphi_v(x,y)$ are functions of $x^2+y^2$, the three requirements are filled.
\end{prop}
The proof is straightforward in the frequency domain, using a stricly analogous proof as in subsection \ref{sec:proof}. Convergence is related to the fact that correction terms of the form \eqref{gains:ff}-\eqref{gains:gg} have  positive Fourier coefficients. Such integral correction terms are not too restrictive, as  convolution with such terms is the integral counterpart of multiplication by a positive scalar gain.  The interest of such observers is that the convergence rate in the frequency domain is speeded up by a factor $\sum_{k=0}^N\lambda_k(p^2+q^2)^k$ on each Fourier coefficient without affecting smoothness of the correction terms. Note that the Proposition remains valid setting  $\varphi_h$ and $\varphi_v$ equal to Dirac functions. Such differential terms can be used in the absence of noise.

\section{Observer design for an oceanography example}\label{oceanography:sec}

The problem considered is the following: the ocean is described by a
simplified shallow-water model. The sea surface height (SSH) is
measured (with noise) everywhere by satellites. The goal is to
estimate the height, and the marine currents (not measured).  There is an increasing need for
such methods in physical oceanography, as the monitoring of the
ocean provides crucial information about climate changes
\cite{Rozier07}, and the
 amount of data available in oceanography has drastically
increased in the last
 years with the use of satellites.

The use of observers for data assimilation in oceanography goes by the name of ``nudging". Indeed the standard nudging algorithm is viewed either as
applying a Newtonian recall of the state value towards its direct
observation \cite{Anthes} or as using observers of the Luenberger, or extended Kalman
filter type for data assimilation \cite{Luenberger,Kalman}. The correction gain is
usually chosen by numerical experimentation. The nudging (i.e. observer) method is known to be
much more economical, computationally  speaking, than variational data
assimilation methods \cite{LeDimet,LeDimetTalagrand}.

Observers of the Kalman filter type are designed to provide, for
each time step, the optimal estimate (i.e. of minimal error
variance) of the system state, by using only the previous estimates
of the state and the last observations \cite{Kalman,Evensen}. In the case of a non-linear physical model the
extended Kalman filter only yields an approximation of the optimal
estimate. As the oceanographic models have become very complex in
the recent years, the high computing cost of the extended Kalman
filter can be prohibitive for data assimilation \cite{Rozier07}. The
nudging techniques are Luenberger gain-scheduled observers and the
expression of the gains requires very few (or no) calculations
\cite{Anthes,Verron-Holland}. Our observer is an improvement of these usual techniques.

In this section we consider a simplified oceanic model. The state of
the ocean is the  SSH, and the horizontal speed of the marine
currents. The choice of the orientation  and the origin of the frame
of $\RR^2$  used to express the horizontal coordinates
$(x,y)\in\RR^2$ is arbitrary: the physical problem is invariant by
rotation and translation.

\subsection{Shallow water model}

The shallow water model  is a basic
model usually considered for simple numerical
experiments in oceanography, meteorology or hydrology \cite{pedlosky}, which represents well enough the dynamics of geophysical flows. The
equations are derived from a vertical integration of the
three-dimensional fields, under the hydrostatic approximation,
i.e. neglecting  the vertical acceleration. We consider here the standard
shallow water model of Jiang et al \cite{jiang}. For deeper water,
this model can be adapted into a multi-layer model, each layer being
described by a shallow water model, with some additional terms
modeling stress and friction due to the other layers.

The fluid is made of a layer of constant density $\rho$ with
varying thickness (or height) $h(x,y,t)$, covering a deeper layer of
density $\rho+\Delta\rho$. The domain is still rectangular: $0\le x\le L$
and $0\le y\le L$, where $x$ and $y$ are the cartesian coordinates
corresponding to East and North respectively.
The equations write:
\begin{align}
&\frac{\partial (hv)}{\partial t}+(\nabla\cdot
(hv)+(hv)\cdot\nabla)v=-g'h\nabla h-\textbf{k}\times
f(hv)\nonumber\\
\label{shallow:eq1}& \qquad\qquad+(A\nabla^2-R)(hv)+\tilde{\tau}\textbf{i}/\rho,\\
\label{shallow:eq2}&\frac{\partial h}{\partial t}=-\nabla\cdot (hv),
\end{align}
where $hv=h(v_x \textbf{i} + v_y \textbf{j})$ is the horizontal
transport, with $\textbf{i}$ and $\textbf{j}$ pointing towards East
and North respectively, $f=f_0+\beta y$ is the Coriolis
parameter (in the $\beta$-plane approximation), $\textbf{k}$ is the
upward unit vector, and $g'$ is the reduced gravity.  The ocean is
driven by a zonal wind stress $\tilde{\tau}\textbf{i}$ modeled as a
body force, and $\tilde{\tau}$ is known. Finally, $R$ and $A$
represent friction and lateral viscosity. No-slip boundary conditions are imposed, i.e. $v=0$ on the boundary of the domain (see paragraph \ref{sec:SV} for more details about the boundary conditions).

We briefly describe the numerical schemes used for the resolution of these equations (as well as the linearized Saint-Venant equations, and all observer equations). We refer to \cite{these-durbiano} for more details. We consider a leap-frog method for time discretization of the equations, controlled by an Asselin time filter \cite{asselin}. The equations are then discretized on an Arakawa C grid \cite{arakawa}, with $N\times N$ points: the velocity components $v_x$ and $v_y$ are defined at the center of the edges, and the height is defined at the center of the grid cells. Then, the vorticity and Bernoulli potential are computed at the nodes and center of the cells respectively. This scheme is known to give stable and accurate results.

We assume that the physical system is observed by several satellites
that provide (noisy) measurements of the SSH $h(x,y,t)$ for all $x,y,t$.
Within the framework of data assimilation for geophysical fluids,
the goal is to estimate all the state variables $v(x,y,t)$ and
$h(x,y,t)$ (velocity of the marine streams, and SSH respectively) at
any point $(x,y)\in [0,L]^2$ of the domain. We finally consider that
all the other parameters are known.

As previously mentioned, if the height is only measured on a discrete set, one usually gathers these sets over a standard time period (e.g. 1 day, or a few days, for oceans), and then interpolates this set in order to have a full observation of the height. With this approach, we can consider that $h$ is observed everywhere in space, but only at some discrete times. The correction term in the observer equations can then be added only at these observation times. Of course, the convergence of the observer towards the real solution is slower than for full observations (in time), and the solution at convergence is less precise, but from the numerical point of view, the method is still applicable.

\subsection{Model symmetries}

The unit vectors $\textbf{i}$ and $\textbf{j}$ are pointing East and
North respectively. This choice is arbitrary, and the equations of
fluid mechanics  are
invariant under the action of $SE(2)$. Considering the transformations \eqref{new:coord:eq1}-\eqref{new:coord:eq2}-\eqref{new:coord:eq3}, the equations in the new
coordinates are unchanged. Indeed letting $\textbf{K}=\textbf{k}$
and $\textbf{I} = R_\theta\textbf{i}$ we have:
\begin{eqnarray}
&&\frac{\partial (HV)}{\partial t}+(\nabla\cdot
(HV)+(HV)\cdot\nabla)V=-g'H\nabla H\nonumber\\
&& \qquad -\textbf{K}\times f(HV)+(A\nabla^2-R)(HV)+\tilde{\tau}\textbf{I}/\rho,\\
&&\frac{\partial H}{\partial t}=-\nabla\cdot (HV).
\end{eqnarray}
where the square domain
$D=[0,L]^2\subset\RR^2$ is replaced by
$\left(R_\theta D+(x_0,y_0)\right)\subset\RR^2$. The boundary conditions are obviously unchanged.

\subsection{Symmetry-preserving nudging}

An observer for the system \eqref{shallow:eq1}-\eqref{shallow:eq2}
 (nudging estimator) systematically writes :
\begin{eqnarray}
&&\frac{\partial (\hat h\hat v)}{\partial t}+(\nabla\cdot (\hat h\hat
v)+(\hat h\hat v)\cdot\nabla)\hat v=-g'\hat h\nabla \hat
h-\textbf{k}\times
f(\hat h\hat v) \nonumber \\
&&\qquad+(A\nabla^2-R)(\hat h\hat v) +\tilde{\tau}\textbf{i}/\rho
+F_{v}(h,\hat v,\hat h), \label{obs1:eq1}\\
&&\frac{\partial \hat h}{\partial t}=-\nabla\cdot (\hat h\hat
v)+F_h(h,\hat v,\hat h), \label{obs1:eq2}
\end{eqnarray}
with $\hat v=0$ on the boundary of the domain, and where the
correction terms vanish when the estimated height $\hat h$ is equal
to the observed height $h$:
 $F_{v}(h,\hat v,h)=0,\ F_{h}(h,\hat v,h)=0.$

As the system possesses the same symmetries as \eqref{saint-venant:eq1}-\eqref{saint-venant:eq2}, we get a large class of $SE(2)$-invariant candidate  correction terms  given by \eqref{fv:eq}-\eqref{fh:eq}. In subsection \ref{num:full:sec}, devoted to numerical experiments, we focus on the particular choice of Section \ref{obs:sec}, i.e., $F_h=\varphi_h * (h-\hat h)$ and $F_v=\varphi_v*\nabla(h-\hat h)$ with the kernels given by  \eqref{gains:eq1}-\eqref{gains:eq2}.  Even if we have no proof of convergence for the observer  \eqref{obs1:eq1}-\eqref{obs1:eq2} with those correction terms, it is clear from the following numerical experiments that the observer competes with  standard oceanography variational methods, and  is remarkably robust to noise.

\section{Numerical simulations}\label{num:sec}

In this section, we report the results of many numerical simulations
on both the linearized and non-linear shallow water models, in order
to illustrate the interest of such symmetry-preserving observers.
 First the theoretical properties of the observer proved in Section \ref{saint-venant:sec} are illustrated by simulations (subsection \ref{num:linearized:sec}). Then we show on the realistic full non-linear
shallow water model of Section \ref{oceanography:sec} that the observer yields better results than the standard nudging techniques (subsection \ref{num:full:sec}).

\subsection{Linearized simplified system}\label{num:linearized:sec}

We first consider a non-linear shallow water model, in a quasi-linear situation
(small velocities, and height close to the equilibrium height) given
by equations \eqref{saint-venant:eq1}-\eqref{saint-venant:eq2}.
%
The corresponding observer is solution of equations
\eqref{observateur:eq1}-\eqref{observateur:eq2}.
\paragraph*{Remark} Note that in the degenerate case where $\phi_h=K_h\delta_0$ and
$\phi_v=K_{v}\delta_0$ ($K_h$ and $K_{v}$ are positive scalars), we
find the standard nudging terms \cite{npg}:
\begin{eqnarray}
\frac{\partial \hat h}{\partial t}&=&-\nabla\cdot(\hat h\hat v)+K_h(h-\hat h),\\
\frac{\partial \hat v}{\partial t}&=&-(\hat v\cdot\nabla) \hat v-g\nabla
\hat h+K_{v} \nabla(h-\hat h).
\end{eqnarray}

\subsubsection{Model parameters} \label{parameter}

The numerical experiments are performed on a square box, of
dimension $2000$ km\,$\times 2000$ km. The equilibrium height is
$\bar{h}=500$ $m$, and the equilibrium longitudinal and transversal
velocities are $\bar{v}_x=\bar{v}_y=0\ m.s^{-1}$. We consider a
regular spatial discretization with $81\times 81$ grid points. The
corresponding space step is $25$ km. The time step is half an hour
(1800 seconds), and we have considered time periods of $1$ to $4$
months ($1440$ to $5760$ time steps).

The reduced gravity is $g=0.02\ m.s^{-2}$. The height varies between
497.7 and 501.9 m and the norm of the transversal velocity is within
the interval  $\pm 0.008\ m.s^{-1}$. The approximations of the subsection \ref{cv:sec} are valid since $v\ll \sqrt{gh}=3\ m.s^{-1}$ and
$\delta h \ll 500$. The variations of the height and velocities are
indeed of the order of $2$ meters and $0.01\, m.s^{-1}$
respectively. This kind of linearized system with the typical values
above is often considered in geophysical applications, under the
tangent linear approximation, for the estimation of an increment
(instead of the solution itself) \cite{Bennett}.

Concerning the tuning of the gains, we have considered the
convolution kernels defined by equations
\eqref{gains:eq1}-\eqref{gains:eq2}. Recall that $\alpha_h^{-2}$ and
$\alpha_v^{-2}$ represent the characteristic size of the Gaussian
kernel. We will always take $\alpha_h^{-2}=\alpha_v^{-2}=\alpha$. In
most of the experiments below we have $\alpha=1~m^{-2}$.
Unfortunately the weights $\beta_h$ and $\beta_v$ cannot be chosen
too large for numerical reasons, in order to avoid stability issues.
So we always take $\beta_h\leq 10^{-6}$.  Recall that heuristically
the error equation can be approximated by the damped wave equation
\eqref{wav2:eq} with $\bar h\beta_v={L_0^2\omega_0^2-g\bar h}$ and
$\beta_h=2\xi_0\omega_0$.  The weights $\beta_h$ and $\beta_v$ have
two different units, and physical meaning, and (a priori)  there is
no physical reason why they should have approximately the same
magnitude. Nevertheless, for the numerical values of $\beta_h$
considered in this paper, one can check that any value
$0\leq\beta_v\leq \beta_h$ yields a fundamental frequency for the
error system $\omega_0\sqrt{1-\xi_0^2}$ which is close to  the
natural frequency $\sqrt{g\bar h}/L_0$ of the physical system
\eqref{saint-venant:eq1}-\eqref{saint-venant:eq2}. From now on we
will systematically set $\beta_v = 0.1\ \beta_h$, which is
acceptable from a physical point of view, also ensures the
convergence of the observer, and is the largest value of $\beta_v$
which yields numerical stability. Finally, a truncated convolution
integral is used as an approximation of the complete convolution
over the whole domain. The truncation radius is set equal to $10$
pixels in our experiments (further than $10$ pixels away from its
center the Gaussian can be viewed as numerical noise). Close to the boundaries, the convolution integrals are also truncated so that they only cover the domain.


We consider two criteria for quantifying the quality of the
estimation process: the convergence rate of the estimation error,
and the estimation error when convergence is reached. The
initialization of the observer is always
$$
\hat{h}=\bar{h}\ (=500), \quad \hat{v}=\bar{v}\ (=0).
$$
In all the following results, the estimation error is the relative
difference between the true solution ($h$ and $v$) and the observer
solution ($\hat{h}$ and $\hat{v}$):
\begin{equation}
\hspace*{-0.25cm}e_h = \frac{\| (\hat{h}-\bar{h})-(h-\bar{h}) \|}{\| h-\bar{h} \|},
\  e_v = \frac{\| (\hat{v}-\bar{v})-(v-\bar{v}) \|}{\| v-\bar{v}
\|}
\end{equation}
where $\|\,.\,\|$ is the standard $L^2$ norm on the considered
domain. With the previously defined initialization of the observer,
the estimation error at initial time is $e_h(0)=e_v(0)=1$,
corresponding to a 100\% error on the initial conditions. If we
assume that the decrease rate is nearly constant in time, then the
time evolution of the estimation error is given by:
\begin{equation}\label{cv_rate}
e_h(t) = e_h(0) \exp(-c_ht), \ \ e_v(t) = e_v(0) \exp(-c_vt),
\end{equation}
where $c_h$ and $c_v$ are the corresponding convergence rates. In
all the numerical experiments that we have considered, the choice of
the weighting coefficients $\beta_h$ and $\beta_v$ does not modify
the residual estimation errors at convergence. We also noticed in the numerical simulations that the
convergence rates are linearly proportional to $\beta_h$ (and to
$\beta_v=0.1\beta_h$), provided it is not too large.  This is
explained by formula \eqref{eq:upq0} as the Fourier coefficients
$g_{pq}^2$ depend linearly on $\beta_h$.

\subsubsection{Perfect observations}

\begin{figure}
\begin{center}
\includegraphics[angle=270,width=9cm]{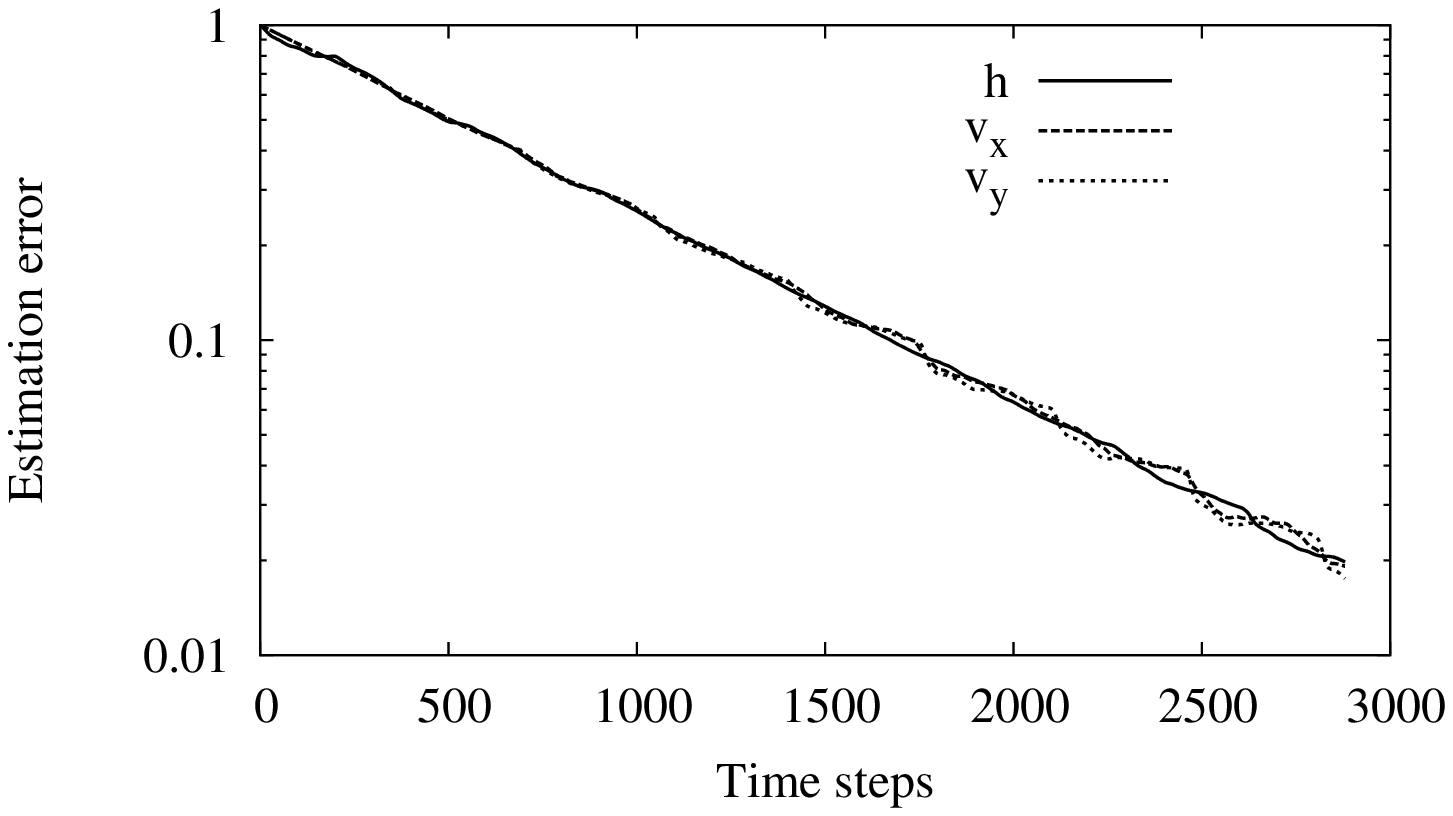}
\caption{Evolution of the estimation error in relative norm versus
the number of time steps, in the case of perfect observations, with
$\alpha_h=\alpha_v=1~m^{-2}$ and $\beta_h=5.10^{-7}~s^{-1}$, and
with a 100\% error on the initial conditions, for the height $h$, longitudinal velocity $v_x$ and transversal
velocity $v_y$.} \label{fig:SV1}
\end{center}
\end{figure}

We first assume that the observations are perfect, i.e. without any
noise. Figure \ref{fig:SV1} shows the estimation error (in relative
norm) versus time (number of time steps), for the three variables:
height $h$, longitudinal velocity $v_x$ and transversal velocity
$v_y$. The kernel coefficients are the following:
$\beta_h = 5.10^{-7}~s^{-1}, \ \beta_v = 0.1 \beta_h =
5.10^{-8}~m.s^{-2}, \ \alpha_h=\alpha_v = 1~m^{-2}.$

This figure shows that the convergence speed is nearly constant in
time, and equation $(\ref{cv_rate})$ is then valid. We can also
deduce the corresponding convergence rates:
$$c_h = 7.57\times 10^{-7}, \quad c_{v_x} = 7.63\times 10^{-7}, \quad c_{v_y}=7.80\times 10^{-7}.$$

In the case of discrete observations, as previously mentioned, we can assume that the height is available everywhere, but not at every time. If for instance we add the correction term to the observer equations only every $12$ time steps, the evolution of the estimation error is similar to what is shown on figure \ref{fig:SV1}, with a smaller convergence rate. In this case, the relative error after 2880 time steps is $0.715$ (to be compared with $0.0198$ in the previous situation), and the convergence rate is approximately $6.47\times 10^{-8}$, which is $11.7$ times smaller than the convergence rate in the full observation case. We could have expected a ratio of $12$, as the corrections are applied only every $12$ time steps. We can conclude that from the numerical point of view, discrete observations in time do not degrade the method.

From an application viewpoint, it is interesting to see that the velocity $v$ is also corrected with a comparable
convergence rate, as predicted by the theory above. Even if it is standard in automatic control theory, in most data assimilation
processes only a few variables of the system are observed
\cite{Anthes,Verron-Holland,npg}. We showed (at least in the
linear case) that  all the variables are observable indeed.

The estimation error at convergence has
the following values:
$$e_h = 7.92\times 10^{-8}, \quad e_{v_x} = 2.11\times 10^{-4}, \quad e_{v_y} = 4.71 \times10^{-5}.$$
From a theoretical point of view, it should converge to $0$. Several
reasons explain this difference with the theory. The numerical
non-linear system considered is not exactly described by its
first-order approximation. Moreover the numerical schemes and
numerical noise do not allow the observer solution to reach exactly
the observed trajectory. Note that the small oscillations in the
decrease of the estimation error can be explained by the oscillatory
behavior described by \eqref{eq:upq0}. Numerically speaking, the
fact that the model has nearly no diffusion (no theoretical
diffusion, and almost no numerical diffusion) can also contribute to
this oscillatory phenomena.

Finally, we compare our observer to the standard nudging algorithm,
by choosing a large value for $\alpha_h$ and $\alpha_v$. Numerically
we have set
$$\alpha_h = \alpha_v = 1000~m^{-2}.$$
The decrease rate and estimation error at convergence are summarized
in table \ref{tab:SV} along with the previous results. The decrease
rate of our observer is $2.7$ to $3$ times bigger. But assuming the
solution $(h,v)$ is constant (which is nearly true), the convolution
with a Gaussian kernel of size $1$ or with a Dirac produces the same
effect, with a $\pi$ factor (as $\int_{\mathbb{R}^2} e^{-(x^2+y^2)}
dx\,dy = \pi$). Numerically, the factor is a little bit smaller, as
the solution is not constant. We also see that the estimation error
at convergence is a little larger for  $\alpha$ large, probably because some
numerical noise  is smoothed by the convolution.
\begin{table}
\begin{center}
\begin{tabular}{l|c|c}
Size of the & Decrease rate & Estimation error at convergence \\
Gaussian kernel & ($h$, $v_x$, $v_y$) & ($h$, $v_x$, $v_y$) \\
\hline
& $7.58\times10^{-7}$ & $7.92\times10^{-8}$ \\
$\alpha_h = \alpha_v = 1$ & $7.63\times10^{-7}$ & $2.11\times10^{-4}$ \\
& $7.80\times10^{-7}$ & $4.71\times10^{-5}$ \\
\hline
& $2.49\times10^{-7}$ & $1.02\times10^{-7}$ \\
$\alpha_h = \alpha_v = 10^3$ & $2.61\times10^{-7}$ & $2.65\times10^{-4}$ \\
& $2.87\times10^{-7}$ & $6.12\times10^{-5}$
\end{tabular}
\caption{Decrease rate and value at convergence of the estimation
error, for the three variables $h$, $v_x$ and $v_y$, for two
different sizes of the Gaussian kernel, in the case of perfect observations.} \label{tab:SV}
\end{center}
\end{table}
\subsubsection{Noisy observations}

\begin{figure}
\begin{center}
\includegraphics[angle=270,width=9cm]{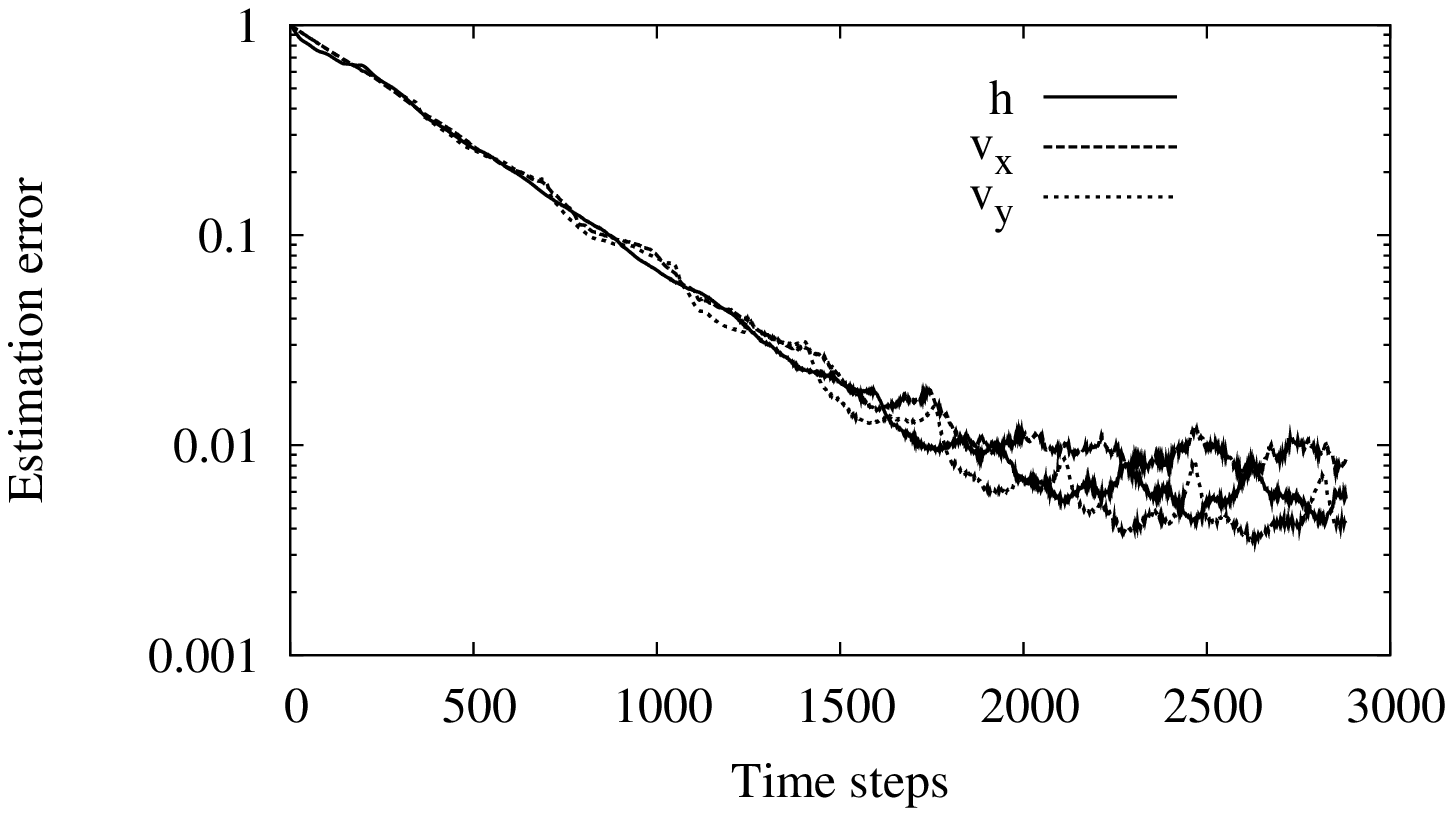}
\caption{Evolution of the estimation error in relative norm versus
the number of time steps, in the case of noisy observations (20\%
noise), with $\alpha_h=\alpha_v=1~m^{-2}$ and $\beta_h =
2.10^{-7}~s^{-1}$, for the three variables: height $h$, longitudinal
velocity $v_x$ and transversal velocity $v_y$.} \label{fig:SV3}
\end{center}
\end{figure}

We now assume that the height $h$ cannot be observed properly, and
instead of $h$, we observe $h+\varepsilon$ where $\varepsilon$
represents the observation noise on $h$. We assume that
$\varepsilon$ is Gaussian with zero mean (white noise is standard in
oceanography \cite{Evensen}), and a standard deviation of $20$ to
$40\%$ of the standard deviation of the height $h$. Thus a $0.2$
relative estimation error means that the estimated value $\hat h$ is closer
to the true height $h$ than to the observed height $h+\varepsilon$.
Figure \ref{fig:SV3} shows similar experiments as previously
described, in the case of noisy observations, for
$\beta_h=2.10^{-7}~s^{-1}$ and $\alpha=1~m^{-2}$. The global
behaviour of the solution is unchanged (constant decrease until
stabilization). The decrease rate and value at convergence of the
estimation error for $\alpha=0.5$, $1$ and $10^3~m^{-2}$ are
summarized in table \ref{tab:SV_bruit}.

\begin{table}
\begin{center}
\begin{tabular}{l|c|c}
Size of the & Decrease rate & Estimation error at convergence \\
Gaussian kernel & ($h$, $v_x$, $v_y$) & ($h$, $v_x$, $v_y$) \\
\hline
& $1.49\times10^{-6}$ & $4.43\times10^{-3}$ \\
$\alpha_h = \alpha_v = 0.5$ & $1.40\times10^{-6}$ & $7.51\times10^{-3}$ \\
& $1.42\times10^{-6}$ & $4.06\times10^{-3}$ \\
\hline
& $7.55\times10^{-7}$ & $5.92\times10^{-3}$ \\
$\alpha_h = \alpha_v = 1$ & $7.44\times10^{-7}$ & $1.04\times10^{-2}$ \\
& $7.44\times10^{-7}$ & $5.53\times10^{-3}$ \\
\hline
& $2.45\times10^{-7}$ & $1.70\times10^{-2}$ \\
$\alpha_h = \alpha_v = 10^3$ & $2.49\times10^{-7}$ & $3.02\times10^{-2}$ \\
& $2.48\times10^{-7}$ & $1.59\times10^{-2}$
\end{tabular}
\caption{Decrease rate and value at convergence of the estimation
error, for the three variables $h$, $v_x$ and $v_y$, for three
different sizes of the Gaussian kernel, in the case of noisy
observations (20\% noise).} \label{tab:SV_bruit}
\end{center}
\end{table}

There is still a ratio of nearly $\pi$ between the decrease rate for
$\alpha$ large and $\alpha=1~m^{-2}$. $\alpha=0.5~m^{-2}$ seems to
be an optimal value for the parameter $\alpha$: it is large enough
to smooth efficiently the noise, and we checked that the decrease
rate is not much larger when we take smaller values of $\alpha$.
Thus we see it is useless to correct the estimation at one point
with values of $h$ which are too far away from this point. In
comparison with the case of perfect observations, the decrease rate
is remarkably unaffected by the presence of noise.

The estimation error at convergence is much larger than in the case
of perfect observations. Nevertheless,  all variables  have been
identified with less than $1\%$ of error. We see the interest of the
convolution as the error at convergence is $3$ to $4$ times smaller
with $\alpha\approx 1$ than with $\alpha=1000$. This is due to the
fact that the term $\nabla(\hat h-h)$ is very noisy when it is not
directly filtered, as it is the case in the standard nudging
algorithm (or extended Kalman filter).

\subsection{Full nonlinear shallow water model}\label{num:full:sec}

We now consider the full shallow water model, with the Coriolis
force, friction, lateral viscosity, and wind stress (see equations
\eqref{shallow:eq1}-\eqref{shallow:eq2}). We also consider large
velocities and height variations, with still the same equilibrium
point: $\bar{h}=500$, $\bar{v}_x=\bar{v}_y=0$. The size of the
domain and the time and space steps remain the same as in the
previous experiments (see section \ref{parameter}), the other
physical parameters being:
\begin{align*}
&f_0=7.10^{-5} s^{-1}, \quad \beta=2.10^{-11} m^{-1}.s^{-1}, \quad R=9.10^{-8},\\
&A=5\, m^2.s^{-1}, \quad \tilde{\tau}_{max}= 0.05\, s^{-2}.
\end{align*}
The nonlinear observer is given by equations
\eqref{obs1:eq1}-\eqref{obs1:eq2}, with $F_h=\varphi_h * (h-\hat h)$ and $F_v=\varphi_v*\nabla(h-\hat h)$, where
$\varphi_h$ and $\varphi_v$ correspond to
\eqref{gains:eq1}-\eqref{gains:eq2}. It is shown in the appendix  that this model reproduces
quite well the evolution of a fluid in the northern hemisphere.
\subsubsection{Perfect observations}

In order to make the paper not too long, we do not provide the figures and tables corresponding to the case of perfect observations.
We consider the same convolution kernels as in the experiments on
the approximated system above, with the same reference parameters
$\beta_h = 5.10^{-7}~s^{-1}$ and $\beta_v = 0.1 \beta_h$. Many curves showing the estimation error versus
time, for the three variables $h,v_x,v_y$, have been
obtained with several values of $\alpha$. The convergence speeds for
$h,v$ are always constant only at the beginning, and decrease continuously
to $0$ after the error goes under some threshold.

Simulations showed that the final estimation error is
much larger than in the previous experiments. Nevertheless, for $\alpha_h=\alpha_v=1 ~m^{-2}$ the height estimation error is close to
$1\%$, which is a very good result, considering the high turbulence

of the model. The velocity is partially identified (with $12$ to
$15\%$ of error in the best situations). The convergence rates are a
little bit larger than in the linearized case  (around $1.10^{-6}$ for $\alpha_h=\alpha_v=1 ~m^{-2}$). The behaviour between
the standard Gaussian convolution ($\alpha=1~m^{-2}$) and the Dirac
convolution ($\alpha=10^3~m^{-2}$) is comparable to the previous
experiments.

\subsubsection{Noisy observations}

\begin{figure}
\begin{center}
\includegraphics[angle=270,width=9cm]{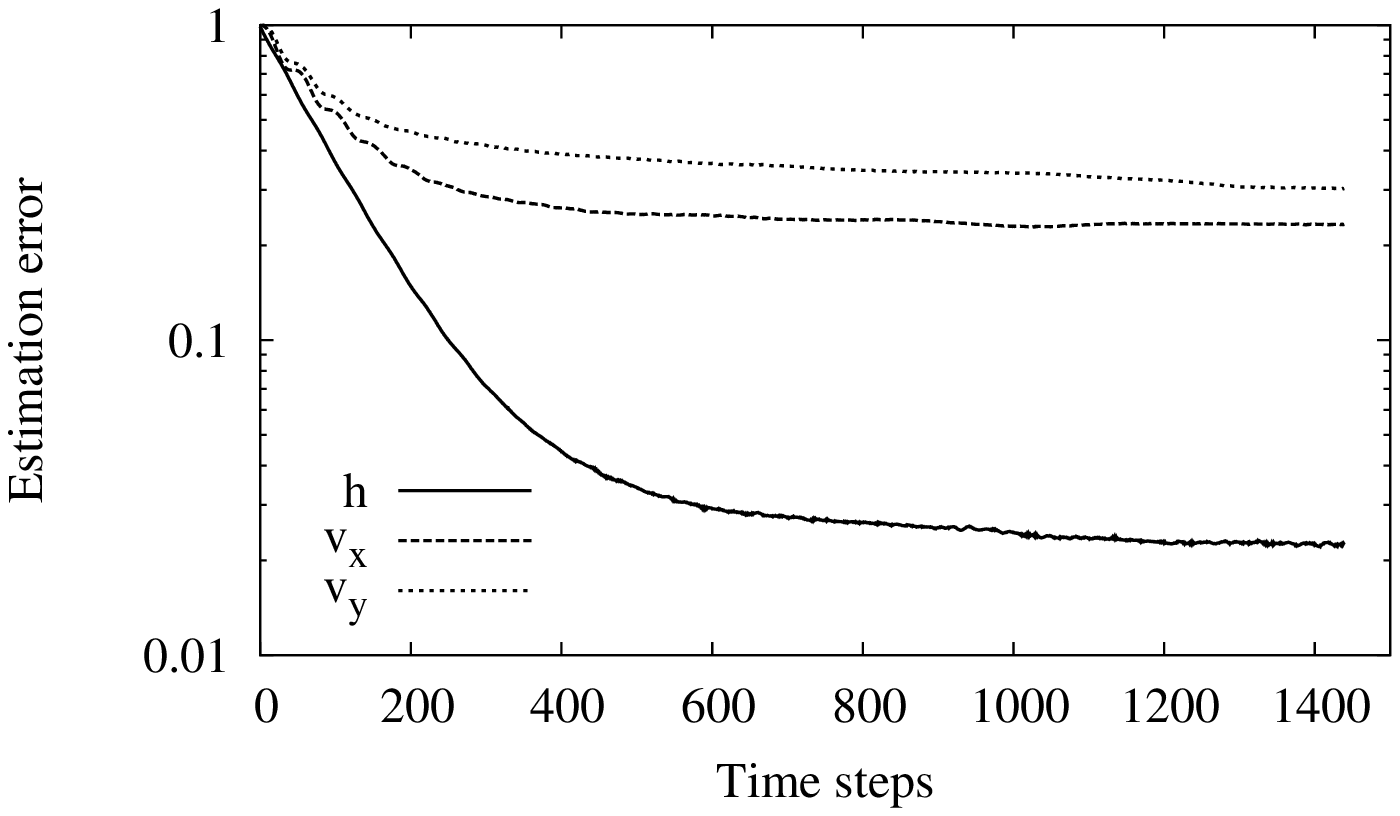}
\caption{Full non-linear model: evolution of the estimation error in relative norm versus
the number of time steps, in the case of noisy observations (20\%
noise), for $\beta_h = 5.10^{-6}~s^{-1}$ and
$\alpha_h=\alpha_v=10^3~m^{-2}$, for the three variables: height
$h$, longitudinal velocity $v_x$ and transversal velocity $v_y$.}
\label{fig:SW2}
\end{center}
\end{figure}

The results are given by figure \ref{fig:SW2} and table \ref{tab:SW_bruit}. As in the linearized situation,  $h+\varepsilon$ is measured, where
$\varepsilon$ is assumed to be white. In our experiments, the
standard deviation of $\varepsilon$ is nearly $20\%$ of the standard
deviation of $h$ (around the equilibrium state $\bar{h}=500$).

\begin{table}
\begin{center}
\begin{tabular}{l|c|c}
Size of the & Decrease rate & Estimation error at convergence \\
Gaussian kernel & ($h$, $v_x$, $v_y$) & ($h$, $v_x$, $v_y$) \\
\hline
& $2.74\times10^{-6}$ & $1.71\times10^{-2}$ \\
$\alpha_h = \alpha_v = 0.5$ & $1.87\times10^{-6}$ & $1.72\times10^{-1}$ \\
& $1.62\times10^{-6}$ & $2.21\times10^{-1}$ \\
\hline
& $1.36\times10^{-6}$ & $1.57\times10^{-2}$ \\
$\alpha_h = \alpha_v = 1$ & $9.65\times10^{-7}$ & $1.30\times10^{-1}$ \\
& $8.38\times10^{-7}$ & $1.59\times10^{-1}$ \\
\hline
& $4.42\times10^{-7}$ & $2.26\times10^{-2}$ \\
$\alpha_h = \alpha_v = 10^3$ & $2.98\times10^{-7}$ & $2.25\times10^{-1}$ \\
& $2.55\times10^{-7}$ & $3.04\times10^{-1}$
\end{tabular}
\caption{Full non-linear model: decrease rate and value at convergence of the estimation
error, for the three variables $h$, $v_x$ and $v_y$, in the case of noisy
observations (20\% noise).} \label{tab:SW_bruit}
\end{center}
\end{table}

The estimation error in the case of noisy observations is nearly
$1.5$ times larger than for perfect observations, both for
 $\alpha=10^3~m^{-2}$ and $\alpha=1~m^{-2}$. The observer has a
relative insensitivity with respect to the presence
of observation noise, as the level of noise is $20\%$, and the
estimation errors are nearly $2\%$ for $h$ and $13$ to $30\%$ for
the velocity. In this case, the best results have been obtained for
$\alpha=1~m^{-2}$, improving the results of the nudging algorithm
($\alpha=10^3~m^{-2}$) of $33$ to $50\%$. These results clearly show
the interest of a Gaussian kernel applied to the correction term, in
order to smooth the noisy observations (or the numerical noise).

The estimation error is of the order of $15\%$ for the velocity at convergence. For instance, if we compare with the standard variational algorithm 4D-VAR \cite{LeDimetTalagrand}, in this kind of situation with noisy observations, the relative error of the velocity at convergence is a little bit larger for 4D-VAR, approximately $18\%$ to $20\%$. Although the results are of the same order, the computing time is totally different: the 4D-VAR needs a few tens of iterations, each iteration consisting of one resolution of the direct model and one resolution of the adjoint model over the time period. Thus the 4D-VAR needs much more computing time than our observer for similar results.

\begin{figure*}
\begin{center}
\hspace*{-0.3cm}\includegraphics[width=\textwidth]{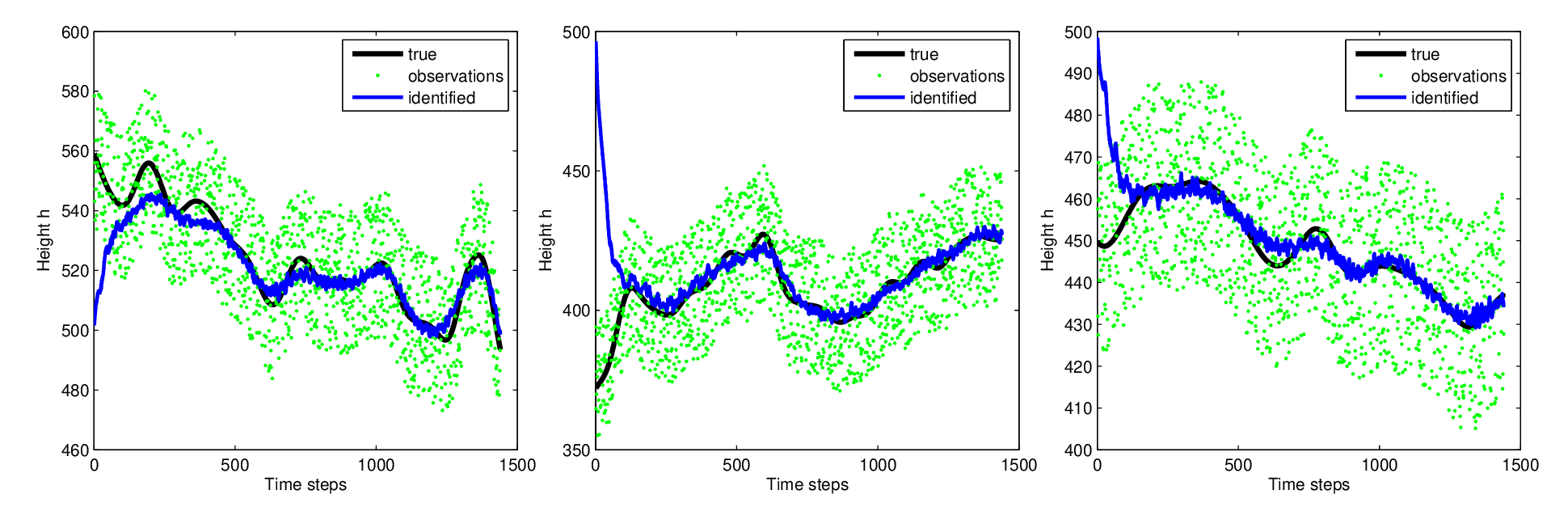}\\[-0.5cm]
\caption{Evolution of the true height, the observed (noisy) height, and the identified (observer) height versus time, for three different points of the domain, located along the energetic current in the middle of the domain.}
\label{fig:SWdiff}
\end{center}
\end{figure*}

In order to show how the observer converges towards the true height, we show on figure \ref{fig:SWdiff} the real height, the observed (noisy) height, and our observer as a function of time, for three different locations inside the domain. These three points are approximately located along the energetic current in the middle of the domain (see figures \ref{fig:SW3} and \ref{fig:SW4}). We can see that after 100 to 400 time steps, the observer is very close to the true height. We can also see that the observation noise is almost totally filtered. In the case of perfect observations (without any noise), the convergence towards the true height is also achieved after a few hundreds of time steps, and the identified height has fewer oscillations around the true height.

\section{Conclusion}\label{ccl:sec}

In this paper, we have defined a class of symmetry-preserving non-linear
observers for a simplified shallow water model. We
proved the asymptotic convergence to zero of the state-error around a steady-state. Many
numerical simulations show the interest of such a choice of
invariant gains. This paper gives insight in the field of non-linear
observers for infinite dimensional systems, where few methods
are available.
The observer provides better results than the nudging (Luenberger observer), even
on the nonlinear system, as the error converges faster, the residual
error is smaller, and the observer is much more robust to noise. The
 correction terms used in this paper are based on integrals over space, and filter the noise better than those
of the usual extended Kalman filter-type estimators. Our observer has several advantages compared to EKF. First the computational cost is much
smaller (as long as the Gaussian kernel is set equal to zero
wherever its value is negligible, see Section \ref{num:sec}). This
is important as in infinite dimensional systems, the computational
cost of the Kalman filter can be prohibitive, as well as the cost of
optimal techniques (especially in oceanography \cite{LeDimet}). In particular the observer was compared to the standard variational method 4D-Var, and the computing time is much smaller.  Moreover the tuning of
the gains of our observer is very easy as it depends on a very
reduced number of parameters which have a physical meaning (thus the observer is much easier to implement). It is
precisely the use of the physical structure of the system which
allows us to reduce the degrees of freedom in the gain design.
Finally, to the author's knowledge, there is no proof of convergence
of the Kalman filter for infinite dimensional non-linear systems.
Note that we also showed, both on theoretical and numerical points
of view, that the non-observed variables can be corrected, which is still a challenge in geophysics \cite{Bennett}.

We have the following additional comments:\begin{enumerate}
\item Another direction for
future work would be to make numerical experiments  on back and
forth nudging based on our observer. The observer can  easily be
adapted in reverse time  with $\varphi_h\mapsto -\varphi_h$
and $\varphi_v$ unchanged. This new observer-based method has recently appeared, see e.g. \cite{npg} for more details.
\item In this paper we mostly
considered time and space continuous measurements. Some other experiments
could be carried out in the case of sparse observations, both in
time and space.
\end{enumerate}
As a more general concluding remark, although this
paper is only concerned with examples, it yields a systematical way
to take advantage of the rotational invariance   of the Laplacian,
and provides a method for the convergence analysis. A large class of sensible observers can be derived from a list of three formal requirements of subsection \ref{class:sec}. This technique can be an interesting guideline to derive
 novel non-linear observers for
 other estimation problems from physics and
engineering, where the models are based on PDEs (wave equation,
 heat equation) and possess symmetries.

\section*{Acknowledgments}
This paper is supported by the CNRS. The authors would like to
thank  Pierre Rouchon and all members of the LEFE BFN project for useful
discussions.

\bibliographystyle{plain}

\appendix\label{app:sec}

\begin{figure}
\begin{center}
\hspace*{-0.3cm}\includegraphics[width=9.5cm]{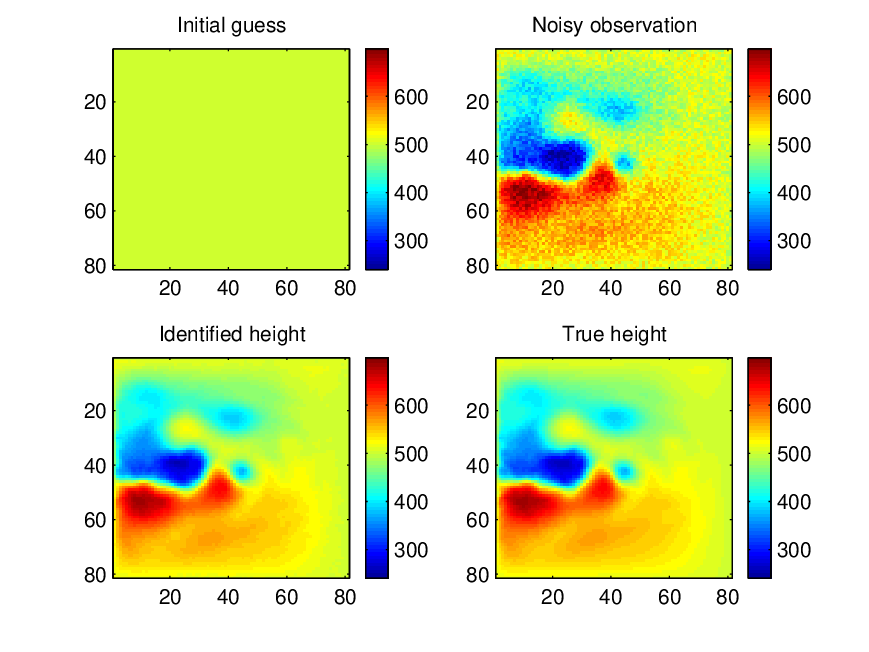}\\[-0.4cm]
\caption{Identification process for the height, in meters: initial
guess ($\hat h(0)=\bar h$); noisy observation at final time
($h(T)+\varepsilon$, with $T=1440$ time steps); identified height at
final time ($\hat h(T)$); true height at final time ($h(T)$).}
\label{fig:SW3}
\end{center}
\end{figure}

\begin{figure}
\begin{center}
\hspace*{-0.3cm}\includegraphics[width=9.5cm]{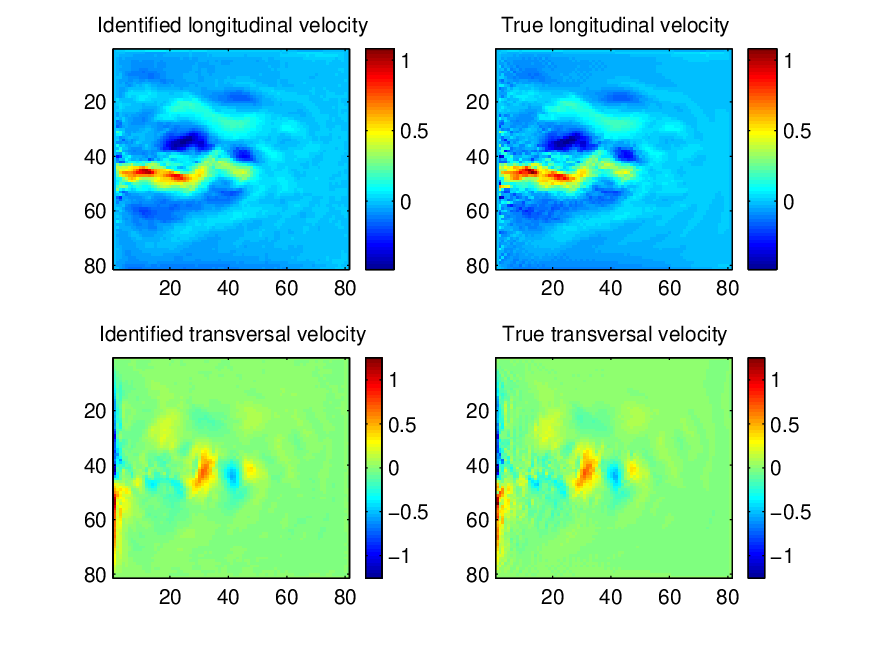}\\[-0.4cm]
\caption{Identification process for the velocity, in $m.s^{-1}$: identified longitudinal (resp. transversal) velocity at final time ($\hat v(T)$); true longitudinal (resp. transversal) velocity at final time ($v(T)$).} \label{fig:SW4}
\end{center}
\end{figure}

In this section we show that the model considered in this paper reproduces
quite well the evolution of a fluid in the northern hemisphere
(e.g. Gulf Stream, in the case of the North Atlantic ocean), with realistic velocities and
dimensions \cite{pedlosky}, and that the observer identifies very well the main currents. Figures \ref{fig:SW3} and \ref{fig:SW4} illustrate the identification process for both the height and velocity in the case of noisy observations, for $\alpha_h=\alpha_v=1$ (second case of table \ref{tab:SW_bruit}). We do not use any a priori information, as the initial guess is $\hat h = \bar h = 500$ meters (top left image of figure \ref{fig:SW3}), and $\hat v = \bar v = 0\ m.s^{-1}$. Figure \ref{fig:SW3} shows on the top right the noisy observation $h+\varepsilon$ of the height at the final time $T=1440$ time steps. It should be compared to the bottom right image, showing the true height $h$ at the same time.
The difference between these two images corresponds to the white Gaussian noise $\varepsilon$. Finally, the identified height (i.e. the observer $\hat h$ at final time $T$) is shown on the bottom left image of figure \ref{fig:SW3}. These images confirm both the very good identification of the height (as previously seen in table \ref{tab:SW_bruit}) and the noise removal.

Figure \ref{fig:SW4} shows the identified and real components of the
velocity. Note that  $\hat v$ is very close to the real velocity
$v$ at time $T$. This is usually not the case in standard nudging
techniques, where only observed variables are corrected and the identification is based on the model coupling
\cite{Anthes,Verron-Holland,npg,Kalnay}. The main current (corresponding to
the Gulf Stream, in the case of the North Atlantic ocean) is very
well identified. This corresponds to a real need, as in operational
geophysical applications, there are also almost no observations of
the fluid velocity, although it has to be precisely identified
\cite{Bennett}. From table \ref{tab:SW_bruit}, we have previously
seen that the error on the velocity is nearly $15\%$ in this case,
which is quite high. But the main currents are very well identified,
and this is a key-point for improving the quality of the forecasts.

\begin{thebibliography}{10}

\bibitem{adcroft}
A.~Adcroft and D.~Marshall.
\newblock How slippery are piecewise-constant coastlines in numerical ocean models?
\newblock {\em Tellus} 50(1):95--108, 1998.

\bibitem{Aghannan}
N.~Aghannan and P.~Rouchon.
\newblock On invariant asymptotic observers.
\newblock In {\em Proc. of the 41st IEEE Conf. on Decision and Control},
  pages 1479--1484, 2002.

\bibitem{aghannan-rouchon-ieee03}
N.~Aghannan and P.~Rouchon.
\newblock An intrinsic observer for a class of lagrangian systems.
\newblock {\em IEEE Trans. on Automatic Control}, 48(6):936--945, 2003.

\bibitem{lions}
L.~Alvarez, F.~Guichard, P.-L. Lions, and J.-M. Morel.
\newblock Axioms and fundamental equations of image processing.
\newblock {\em Arch. Rational Mech. Anal.}, 123:199--257, 1993.

\bibitem{arakawa}
A.~Arakawa and V.~Lamb.
\newblock Computational design of the basic dynamical processes of the {UCLA} general circulation model.
\newblock  {\em Methods in Computational Physics}, 17:174--267, Academic Press, 1977.

\bibitem{asselin}
R.~Asselin.
\newblock Frequency filter for time integrations.
\newblock {\em Mon. Wea. Rev.}, 100:487--490, 1972.

\bibitem{npg}
D.~Auroux and J.~Blum.
\newblock A nudging-based data assimilation method for oceanographic problems:
  the {B}ack and {F}orth {N}udging ({BFN}) algorithm.
\newblock {\em Nonlin. Proc. Geophys.}, 15:305--319, 2008.

\bibitem{Bennett}
A.~F. Bennett.
\newblock {\em Inverse Modeling of the Ocean and Atmosphere}.
\newblock Cambridge University Press, Cambridge, 2002.

\bibitem{besancon}
G. Besancon, J.F. Dulhoste, D.  Georges.
\newblock Nonlinear observer design for water level control in irrigation canals.
\newblock In {\em 40th IEEE Conf. on Decision and Control}, pp. 4968--4973, 2001.


\bibitem{arxiv-07}
S.~Bonnabel, Ph. Martin, and P.~Rouchon.
\newblock Symmetry-preserving observers.
\newblock {\em IEEE Trans. on Automatic Control}, 53(11):2514--2526, 2008.


\bibitem{arxiv-08}
S.~Bonnabel, Ph. Martin, and P.~Rouchon.
\newblock Non-linear symmetry-preserving observers on Lie groups.
\newblock {\em IEEE Trans. on Automatic Control}, 54(7):1709--1713, 2009.


\bibitem{bonnabel_auto}
S.~Bonnabel, M. Mirrahimi, and P.~Rouchon.
\newblock Observer-based Hamiltonian identification for quantum systems.
\newblock {\em Automatica}, 45: 1144-1155, 2009.


\bibitem{coron-cocv02}
J.~M. Coron.
\newblock {Local controllability of a 1-D tank containing a fluid modeled by
the shallow water equations},
\newblock {\em ESAIM: COCV},
8:513--554, 2002.


\bibitem{coron-et-al-ecc99}
       J.~M. Coron, B. D'Andr\'ea-Novel and G. Bastin.
  \newblock {A Lyapunov approach to control
        irrigation canals modeled by Saint-Venant
        equations},
   \newblock {\em Proc. European Control Conference, Karlsruhe}, 1999.

\bibitem{coron-ieee}
       J.~M. Coron, B. D'Andr\'ea-Novel and G. Bastin.
  \newblock {A Strict Lyapunov Function for Boundary Control of Hyperbolic Systems of Conservation Laws},
   \newblock {\em  IEEE Trans. on Automatic Control} 52(1):2 -- 11, 2007.


\bibitem{deguenon2}
J. Deguenon, G. Sallet and C.~Z. Xu.
\newblock        {A Kalman observer for infinite dimensional skew-symmetric systems with applications to an elastic beam},
\newblock {\em 2nd Int. Symp. Communications, Control, Signal Processing},
2006.

\bibitem{prieur}
V. Dos Santos and C. Prieur.
\newblock        {Boundary control of open channels with numerical and experimental validations},
\newblock {\em IEEE Trans. on Control Syst. Tech.}, 16(6):1252--1264, 2008.


\bibitem{dubois}
F.~Dubois, N.~Petit, and P.~Rouchon.
\newblock Motion planing and nonlinear simulations for a tank containing a
  fluid.
\newblock  {\em European Control Conference, Karlsruhe}, 1999.

\bibitem{these-durbiano}
S.~Durbiano.
\newblock Vecteurs caract\'eristiques de mod\`eles oc\'eaniques pour la r\'eduction d'ordre en assimilation de donn\'ees.
\newblock {\em PhD thesis}, University of Grenoble, 2001.

\bibitem{Evensen}
G.~Evensen.
\newblock The ensemble Kalman filter: Theoretical formulation and practical
  implementation.
\newblock {\em Ocean Dynam.}, 53:343--367, 2003.

\bibitem{GuoShao}
B.~Z.~Guo and Z.~C.~Shao.
\newblock Stabilization of an abstract second order system with application to wave equations under non-collocated control and observations.
\newblock {\em Systems and Control Letters}, 58:334--341, 2009.

\bibitem{Guo}
B.-Z. Guo and C.~Z. Xu.
\newblock The stabilization of a one dimensional wave equation by boundary feedback with noncollocated observation.
\newblock {\em IEEE Trans. on Automatic Control}, 52:371--377, 2007.


\bibitem{Anthes}
J.~Hoke and R.~A. Anthes.
\newblock The initialization of numerical models by a dynamic initialization
  technique.
\newblock {\em Month. Weather Rev.}, 104:1551--1556, 1976.

\bibitem{jiang}
S.~Jiang and M.~Ghil.
\newblock Tracking nonlinear solutions with simulated altimetric data in a
  shallow-water model.
\newblock {\em J. Phys. Oceanogr.}, 27(1):72--95, 1997.

\bibitem{Kalman}
R.~E. Kalman.
\newblock A new approach to linear filtering and prediction problems.
\newblock {\em Trans. ASME - J. Basic Engin.}, 82:35--45, 1960.

\bibitem{Kalnay}
E.~Kalnay.
\newblock {\em Atmospheric modeling, data assimilation and predictability}.
\newblock Cambridge University Press, 2003.

\bibitem{komornik-book-2}
V.~Komornik and P.~Loreti.
\newblock {\em Fourier Series in Control Theory}.
\newblock Springer Monographs in Mathematics. Springer-Verlag, New York, 2005.

\bibitem{lageman}
C. Lageman, R. Mahony, and J. Trumpf.
\newblock Gradient-like observers for invariant dynamics on a Lie group.
\newblock {\em IEEE Trans. on Automatic Control},   55:367--377, 2010.

\bibitem{lasiecka}
I.~Lasiecka and R.~Triggiani.
\newblock {\em Control theory for partial differential equations: continuous and approximation theories}.
\newblock Cambridge University Press, 2000.



\bibitem{LeDimetTalagrand}
F.-X. Le~Dimet and O.~Talagrand.
\newblock Variational algorithms for analysis and assimilation of
  meteorological observations: theoretical aspects.
\newblock {\em Tellus}, 38A:97--110, 1986.

\bibitem{Luenberger}
D.~Luenberger.
\newblock Observers for multivariable systems.
\newblock {\em IEEE Trans. on Automatic Control}, 11:190--197, 1966.

\bibitem{mahony-et-al-dcd05}
R.~Mahony, T.~Hamel, and J.-M. Pflimlin.
\newblock Non-linear complementary filters on the special orthogonal group.
\newblock In {\em  IEEE Trans. on Automatic Control}, 53:1203--1218 ,
  2008.



\bibitem{mareels}
I. Mareels, E. Weyer, S.~K. Ooi, M. Cantoni, Y. Li, and G. Nair.
\newblock Systems engineering for irrigation systems: Successes and challenges.
\newblock In {\em Annual Reviews in Control}, 29(2): 191--204, 2005.


\bibitem{MarSal2007a}
Ph. Martin and E.~Salaun.
\newblock Invariant observers for attitude and heading estimation from low-cost
  inertial and magnetic sensors.
\newblock In {\em 46th IEEE Conf. on Decision and Control}, pages 1039--1045.

\bibitem{olver-book95}
P.~J. Olver.
\newblock {\em Equivalence, Invariants, and Symmetry}.
\newblock Cambridge University Press, 1995.

\bibitem{pedlosky}
J.~Pedlosky.
\newblock {\em Geophysical fluid dynamics}.
\newblock Springer-Verlag, New-York, 1979.

\bibitem{petit-rouchon-ieee02}
N.~Petit and P.~Rouchon.
\newblock Dynamics and solutions to some control problems for water-tank
  systems.
\newblock {\em IEEE Trans. on Automatic Control}, 47(4):594--609, 2002.

\bibitem{Rozier07}
D.~Rozier, F.~Birol, E.~Cosme, P.~Brasseur, J.-M. Brankart, and
J.~Verron.
\newblock A reduced-order kalman filter for data assimilation in physical
  oceanography.
\newblock {\em SIAM Rev.}, 49(3):449--465, 2007.

\bibitem{schwartz-seminaire}
L.~Schwartz.
\newblock Op\'erateurs invariants par rotations. Fonctions m\'etaharmoniques.
\newblock In {\em S\'eminaire Schwartz}, 2(9):1--5, 1954.



\bibitem{krstic}
A. Smyshlyaev and M.~Krstic.
\newblock Backstepping observers for a class of parabolic PDEs.
\newblock {\em Systems Control Letters}, 54: 613--625, 2005.

\bibitem{krstic-05}
R. Vazquez  and M.~Krstic.
\newblock A closed-form observer for the channel flow Navier-Stokes system.
\newblock {\em 44th IEEE Conf. on Decision Control}, 2005.

\bibitem{Verron-Holland}
J.~Verron and W.~R. Holland.
\newblock Impact de donn\'ees d'altim\'etrie satellitaire sur les simulations
  num\'eriques des circulations g\'en\'erales oc\'eaniques aux latitudes
  moyennes.
\newblock {\em Annales Geophysicae}, 7(1):31--46, 1989.



\bibitem{deguenon}
C.~Z. Xu, J. Deguenon and G. Sallet.
\newblock        {Infinite dimensional observer for vibrating systems.}
\newblock {\em 45th IEEE Conf. on Decision Control},
2005.


\bibitem{LeDimet}
X.~Zou, I.~M. Navon, and F.-X. Le~Dimet.
\newblock An optimal nudging data assimilation scheme using parameter
  estimation.
\newblock {\em Q. J. R. Meteorol. Soc.}, 118:1163--1186, 1992.

\end{thebibliography}

\end{document}